\documentclass[12pt]{article}
\usepackage{latexsym, amssymb
}
    \title{{\bf  Open-string vertex algebras, tensor categories and 
operads}}
    \author{Yi-Zhi Huang and Liang Kong}
    \date{}
    \begin{document}
    \bibliographystyle{alpha}
    \maketitle

   \newtheorem{thm}{Theorem}[section]
\newtheorem{defn}[thm]{Definition}
\newtheorem{prop}[thm]{Proposition}
\newtheorem{cor}[thm]{Corollary}
\newtheorem{rema}[thm]{Remark}
\newtheorem{lemma}[thm]{Lemma}
\newtheorem{app}[thm]{Application}
\newtheorem{prob}[thm]{Problem}
\newtheorem{conv}[thm]{Convention}
\newtheorem{conj}[thm]{Conjecture}
\newtheorem{expl}[thm]{Example}
\newcommand{\halmos}{\rule{1ex}{1.4ex}}
\newcommand{\pfbox}{\hspace*{\fill}\mbox{$\halmos$}}

\newcommand{\bnu}{\begin{enumerate}}
\newcommand{\enu}{\end{enumerate}}
 \newcommand{\bea}{\begin{eqnarray}}
\newcommand{\eea}{\end{eqnarray}}
 \newcommand{\nn}{\nonumber \\}
	\newcommand{\nno}{\nonumber}
	\newcommand{\lbar}{\bigg\vert}
\newcommand{\mbar}{\mbox{\large $\vert$}}
	\newcommand{\p}{\partial}
	\newcommand{\dps}{\displaystyle}
	\newcommand{\bra}{\langle}
	\newcommand{\ket}{\rangle}
 \newcommand{\res}{\mbox{\rm Res}}
\renewcommand{\hom}{\mbox{\rm Hom}}
\newcommand{\hol}{\mbox{\rm Hol}}
\newcommand{\dt}{\mbox{\rm Det}}
\newcommand{\edo}{\mbox{\rm End}\;}
 \newcommand{\pf}{{\it Proof.}\hspace{2ex}}
 \newcommand{\epf}{\hspace*{\fill}\mbox{$\halmos$}}
 \newcommand{\epfv}{\hspace*{\fill}\mbox{$\halmos$}\vspace{1em}}
 \newcommand{\epfe}{\hspace{2em}\halmos}
\newcommand{\nord}{\mbox{\scriptsize ${\circ\atop\circ}$}}
\newcommand{\wt}{\mbox{\rm wt}\ }
\newcommand{\swt}{\mbox{\rm {\scriptsize wt}}\ }
\newcommand{\clr}{\mbox{\rm clr}\ }
\newcommand{\id}{\mbox{\rm id}}
\newcommand{\injto}{\hookrightarrow}  
\newcommand{\one}{\mathbf{1}}
\renewcommand{\th}{\theta}
\newcommand{\binom}[2]{{{#1}\choose {#2}}}
\newcommand{\text}[1]{\mbox{\rm #1}}
\newcommand{\substack}[2]{\mbox{\rm {\scriptsize $\begin{array}{c}{ 
#1}\\{ #2}\end{array}$}}}

\newcommand{\C}{\mathbb{C}}
\newcommand{\F}{\mathbb{F}} 
\newcommand{\HH}{\mathbb{H}}
\newcommand{\N}{\mathbb{N}}
\newcommand{\Q}{\mathbb{Q}}
\newcommand{\R}{\mathbb{R}}
\newcommand{\Z}{\mathbb{Z}}
\newcommand{\Y}{\mathcal{Y}}

\begin{abstract}
We introduce notions of open-string vertex algebra, conformal 
open-string vertex algebra and variants of these notions. These are 
``open-string-theoretic,'' ``noncommutative'' generalizations of the 
notions of vertex algebra and of conformal vertex algebra. Given an 
open-string vertex algebra, we show that there exists a vertex algebra, 
which we call the ``meromorphic center,'' inside the original algebra 
such that the original algebra yields a module and also an intertwining 
operator for the meromorphic center. This result gives us a general 
method for constructing open-string vertex algebras. Besides obvious 
examples obtained {}from associative algebras and vertex (super)algebras, 
we give a nontrivial example constructed {}from the minimal model of 
central charge $c=\frac{1}{2}$. We establish an equivalence
between the associative algebras in 
the braided tensor category of modules for a suitable
vertex operator algebra and the grading-restricted conformal 
open-string vertex algebras containing a
vertex operator algebra isomorphic to the given vertex operator algebra.
We also give a geometric and operadic formulation of 
the notion of grading-restricted conformal open-string
vertex algebra, we prove two isomorphism theorems,
and in particular, we show that such an algebra gives
a projective algebra over what we call the 
``Swiss-cheese partial operad.''
\end{abstract}

\renewcommand{\theequation}{\thesection.\arabic{equation}}
\renewcommand{\thethm}{\thesection.\arabic{thm}}
\setcounter{equation}{0}
\setcounter{thm}{0} 
\setcounter{section}{-1}
\section{Introduction}

In the present paper, we introduce and study ``open-string-theoretic,''
``noncommutative'' generalizations of ordinary vertex algebras and
vertex operator algebras, which we call ``open-string vertex algebras''
and ``conformal open-string vertex algebras.'' This is a first step in a
program to establish the fundamental and highly nontrivial
assumptions used by physicists in the study of boundary (or open-closed)
conformal field theories as mathematical theorems and to construct
such theories mathematically. See \cite{H8} and \cite{HK2} for
definitions of open-closed conformal field theory in the spirit of the
definition of closed conformal field theory first given by Segal
\cite{S0}--\cite{S2} and Kontsevich in 1987 and further rigorized by Hu
and Kriz \cite{HK1} recently. More recently, Moore suggested in 
\cite{M3} that in order to
generalize a certain formula relating a nonlinear $\sigma$ model and the
$K$-theory on its target space to conformal field theories without
obvious target space interpretation, one should define some kind of
algebraic $K$-theory for ``open string vertex operator algebras.'' 
We hope that 
the notions and results in the present paper will provide a 
solid foundation for the formulation and study of such a $K$-theory.

Vertex (operator) algebras were introduced in mathematics by Borcherds 
in \cite{B}. They arose naturally in the vertex operator
construction of representations of affine Lie algebras and in the
construction and study of the ``moonshine module'' for the Monster
finite simple group by Frenkel-Lepowsky-Meurman \cite{FLM} and Borcherds
\cite{B}.  The notion of vertex (operator) algebra corresponds
essentially to the notion of what physicists call ``chiral algebra" in
(two-dimensional) conformal field theory,  a fundamental physical
theory studied systematically first by Belavin, Polyakov and
Zamolodchikov \cite{BPZ}. Vertex operator algebras can be viewed as
``closed-string-theoretic'' analogues of both Lie algebras and 
commutative associative algebras, and they play important roles in a
range of areas of mathematics and physics. 

Recently, in addition to the continuing development of (closed)
conformal field theories, boundary conformal field theories (open-closed
conformal field theories) have attracted much attention.  Boundary
conformal field theory was first developed by Cardy in \cite{C0},
\cite{C1} and \cite{C2} and play a fundamental role in many problems in
condensed matter physics. It has also become one of the main tools in
the study of open strings and $D$-branes (certain important
nonperturbative objects in string theory).  Besides the obvious problem
of constructing and classifying open-closed conformal field theories,
the study of $D$-branes in physics and their possible applications in
geometry have led to exciting and interesting mathematical problems. If
open-closed conformal field theories associated to Calabi-Yau manifolds
or other geometric objects are constructed eventually, they will provide
even more powerful tools in geometry than the corresponding closed
conformal field theories (see, for example, the survey \cite{D} by
Douglas).  The paper \cite{M3} by Moore mentioned above gave another
example of the exciting and interesting mathematical problems associated
to open-closed conformal field theories.

In the framework of topological field theories, boundary topological
field theories (open-closed topological field theories) have been
studied in detail by Lazaroiu \cite{L} and by Moore and Segal 
\cite{Mo} \cite{M2}
\cite{S3}. In this topological case, an open-closed topological field
theory is roughly speaking a (typically noncommutative) Frobenius algebra
and a commutative Frobenius algebra equipped with some other data 
and satisfying
suitable conditions. The commutative Frobenius algebra is the state
space for the closed string part of the theory and the (typically
noncommutative) Frobenius algebra is the state space for the open string
part of the theory. 

To construct and study open-closed conformal field theories, one first
has to find the analogues in the conformal case of commutative and
noncommutative associative algebras. Since the corresponding algebras in
the conformal case must be infinite-dimensional, their construction and
study are much more difficult than the topological ones.
In the conformal case, one lesson we have learned {}from various 
methods used by physicists is that the
construction and study of chiral theories are necessary and crucial
steps. If chiral theories are constructed, full theories can be
constructed using unitary bilinear forms on substructures of chiral
theories called ``modular functors.'' In fact, it is also the chiral
theories which are more similar to topological theories than full 
theories. It is
clear that analogues of commutative associative algebras in the
chiral conformal case are vertex (operator) algebras.  To construct and
study open-closed conformal field theories, one first has to answer the
following question: What are the analogues of noncommutative associative
algebras in the conformal case? 

Assuming the existence of the structure of a modular tensor category 
on the category of modules for a vertex operator algebra 
and the existence of conformal blocks with monodromies compatible with the 
modular tensor category, Felder, Fr\"{o}hlich,
Fuchs and Schweigert \cite{FFFS} and Fuchs, Runkel and Schweigert
\cite{FRS1} \cite{FRS2} studied open-closed
conformal field theory using the theory of tensor categories and
three-dimensional topological field theories. They showed the existence
of consistent operator product expansion coefficients 
for boundary and bulk operators.
In particular, special symmetric Frobenius algebras in
the modular tensor categories of modules are proposed
as analogues in the conformal case of (typically
noncommutative)  Frobenius algebras in the topological case. However, since 
these works are based on the fundamental assumptions mentioned above, even 
in the genus-zero case, the
corresponding open-string-theoretic
and noncommutative analogues of vertex operator
algebras have not been fully constructed and studied, and even 
chiral open-closed
conformal field theories on the disks 
(the simplest parts of open-closed conformal field theories)
have not been fully constructed. 

The present paper is a first step in a program for establishing the
fundamental and highly nontrivial assumptions mentioned above as
mathematical theorems, 
using the results on representations of
vertex operator algebras and closed conformal field theories. In
particular, we solve the problem of constructing open-closed conformal
field theories on the disks satisfying certain differentiability and 
meromorphicity conditions by introducing, constructing and studying
open-string vertex algebras, conformal open-string vertex algebras and
some other variants. These algebras are
the open-string-theoretic or noncommutative analogues of
vertex (operator) algebras we are looking for and, as we shall 
discuss in future publications, 
$D$-branes can be formulated and studied as irreducible modules for 
suitable open-string vertex algebras.
Given an open-string vertex algebra, we show
that there exists a vertex algebra, which we call 
the ``meromorphic center,'' inside
the open-string vertex algebra such that the open-string vertex algebra
yields a module and also an intertwining operator for the meromorphic center.
This relation between open-string vertex algebras and the
representation theory of vertex algebras gives us a general method for
constructing open-string vertex algebras.  Besides obvious
examples obtained {}from associative algebras and vertex (super)algebras,
we give a nontrivial one constructed {}from the minimal model of
central charge $c=\frac{1}{2}$. We establish an equivalence between
grading-restricted conformal open-string vertex algebras containing a
suitable vertex operator algebra and associative algebras in the braided
tensor category of modules for the vertex operator algebra. We also
give a geometric and operadic formulation of the notion of
grading-restricted conformal open-string vertex algebra, we prove two
isomorphism theorems (establishing the equivalence of geometric
notions and algebraic notions), and in particular, we
show that such an algebra gives a
projective algebra over what we call the ``Swiss-cheese partial operad."

Here is the organization of the present paper: In Section 1, we
introduce the notions of open-string vertex algebra, conformal
open-string vertex algebra and other variants. The connection between
open-string vertex algebras the representation theory of vertex
(operator) algebras is given in Section 2. Examples of (conformal)
open-string vertex algebras are presented in Section 3. In Section 4, we
show that for a vertex operator algebra satisfying certain finiteness
and complete reductivity properties, associative algebras in the braided
tensor category of modules for the vertex operator algebra are
equivalent to grading-restricted conformal open-string vertex algebras
containing the vertex operator algebra in their meromorphic centers.
The geometric and operadic formulation of the notion of
grading-restricted conformal open-string vertex algebra, the
construction of projective algebras over the Swiss-cheese partial operad
and the proof of the corresponding isomorphism theorems
are given in Section 5.

We shall use $\mathbb{C}$, $\mathbb{H}$, $\overline{\mathbb{H}}$,
$\hat{\mathbb{H}}$ , $\R$, $\R^{\times}$, $\R_{+}$, $\mathbb{Z}$,
$\mathbb{Z}_{+}$ and $\mathbb{N}$ to denote the sets (with structures)
of the complex numbers, the open upper half plane, the closed upper half
plane, the one point compactification of the closed upper half plane,
the nonzero real numbers, the positive real numbers, the integers, the
positive integers and the nonnegative integers, respectively.  For any
$z\in \C^{\times}$ and $n\in \C$, we shall always use $\log z$ and
$z^{n}$ to denote $\log |z|+\arg z$, $0\le \arg z< 2\pi$, and $e^{n\log
z}$, respectively.

\paragraph{Acknowledgment} We would like to thank J\"{u}rgen Fuchs
and Christoph Schweigert for helpful discussions and comments. 
We are also grateful to Jim Lepowsky for comments.
The research of Y.-Z. H.  is supported in
part by NSF grant DMS-0070800. 

\renewcommand{\theequation}{\thesection.\arabic{equation}}
\renewcommand{\thethm}{\thesection.\arabic{thm}}
\setcounter{equation}{0}
\setcounter{thm}{0}

\section{Definitions and basic properties}

We introduce the notion of open-string vertex algebra and its variants
and discuss some basic properties of these algebras in this section.
We assume that the reader is familiar with the basic notions 
and properties in the theory of vertex operator algebras as presented 
in \cite{FLM} and \cite{FHL}.

In the present paper, all vector spaces are over the field $\C$.  For a
vector space $V$, we shall use $V^{-}$ to denote its complex conjugate
space, which is characterized by the fact that if $\sqrt{-1}$ acts as
$J$ on the underlying real vector space of $V$, then $\sqrt{-1}$ acts as
$-J$ on the underlying real vector space of $V^{-}$.  For an $\R$-graded
vector space $V=\coprod_{n\in \R}V_{(n)}$ and any $n\in \R$, we shall
use $P_{n}$ to denote the projection {}from $V$ or
$\overline{V}=\prod_{n\in \R}V_{(n)}$ to $V_{(n)}$.  We give $V$ and its
graded dual $V'=\coprod_{n\in \R}V^{*}_{(n)}$ the topology induced
{}from the pairing between $V$ and $V'$. We also give $\hom(V,
\overline{V})$ the topology induced {}from the linear functionals on
$\hom(V, \overline{V})$ given by $f\mapsto \langle v', f(v)\rangle$ for
$f\in \hom(V, \overline{V})$, $v\in V$ and $v'\in V'$. 

\begin{defn}
{\rm An {\it open-string vertex algebra} 
is an $\R$-graded vector space $V=\coprod_{n\in \R}V_{(n)}$ (graded
by {\it weights}) 
equipped with a {\it vertex map}
\begin{eqnarray*}
Y^{O}: V \times \R_{+} &\to& \hom(V, \overline{V})\nn
(u, r)&\mapsto &Y^{O}(u, r)
\end{eqnarray*}
or equivalently,
\begin{eqnarray*}
Y^{O}: (V\otimes V) \times \R_{+} &\to& \overline{V}\\
(u\otimes v, r)&\mapsto &Y^{O}(u, r)v,
\end{eqnarray*}
a {\it vacuum} $\one\in V$
and an operator $D\in \edo V$ of weight $1$, satisfying 
the following conditions:

\bnu

\item {\it Vertex map weight property}: 
For $n_{1}, n_{2}\in \R$, 
there exist a finite subset $N(n_{1}, n_{2})\subset \R$ 
such that the image of 
$\left(\coprod_{n\in n_{1}+\Z}V_{(n)}\otimes 
\coprod_{n\in n_{2}+\Z}V_{(n)}\right)\times \R_{+}$
under $Y^{O}$ is in 
$\prod_{n\in N(n_{1}, n_{2}) +\Z}
V_{(n)}.$

\item Properties for the vacuum: For any $r\in \R_{+}$,
$Y^{O}(\one,r)=\id_{V}$ (the
{\it identity property}) and  $\lim_{r\rightarrow 0} Y^{O}(u, r)\one$
exists and is equal to $u$
(the {\it creation property}).

\item {\it Local-truncation property for $D'$}: Let 
$D': V'\to V'$ be the adjoint of $D$.
Then for any $v'\in V'$, there exists a positive integer $k$ such that 
$(D')^{k}v'=0$.

\item {\it Convergence properties}: For $v_{1}, \dots, v_{n}, v\in V$
and $v'\in V'$, the series 
\begin{eqnarray*}
\lefteqn{\langle v', Y^{O}(v_{1}, r_1)\cdots Y^{O}(v_{n}, r_n)v\rangle}\nn
&&=\sum_{m_{1}, \dots, m_{n-1}\in \R}
\langle v', Y^{O}(v_{1}, r_1)P_{m_{1}}Y^{O}(v_{2}, r_2)\cdots 
P_{m_{n-1}}Y^{O}(v_{n}, r_n)v\rangle
\end{eqnarray*}
converges absolutely 
when $r_1>\cdots >r_n>0$. For $v_{1}, v_{2}, v\in V$
and $v'\in V'$, the series 
$$\langle v', Y^{O}(Y^{O}(v_{1}, r_0)v_{2}, r_2)v\rangle$$
converges absolutely when $r_{2}>r_{0}>0$. 

\item {\it Associativity}: For $v_{1}, v_{2}, v\in V$
and $v'\in V'$, 
$$\langle v', Y^{O}(v_{1}, r_1)Y^{O}(v_{2}, r_2)v\rangle =
\langle v', Y^{O}(Y^{O}(v_{1}, r_1-r_2)v_{2},r_2)v\rangle$$
for  $r_{1}, r_{2}\in \R$ satisfying 
$r_1>r_2>r_1-r_2>0$.

\item {\it $\mathbf{d}$-bracket property}: Let $\mathbf{d}$ be the grading 
operator on $V$, that is, $\mathbf{d}u=mu$ for $m\in \R$ and
$u\in V_{(m)}$. For $u\in V$ and $r\in \R_{+}$,
\begin{equation} \label{d-bra}
[\mathbf{d}, Y^{O}(u,r)]= Y^{O}(\mathbf{d}u,r)+r\frac{d}{dr}Y^{O}(u,r).
\end{equation}

\item {\it $D$-derivative property}: We still use $D$ to denote
the natural extension of $D$ to $\hom(\overline{V}, \overline{V})$.
For $u\in V$, 
$Y^{O}(u, r)$ as a map {}from $\R_{+}$ to 
$\hom(V, \overline{V})$
is differentiable and 
\begin{equation}\label{D-der}
\frac{d}{dr}Y^{O}(u,r)=[D, Y^{O}(u,r)]
= Y^{O}(Du,r).
\end{equation}

\enu  

{\it Homomorphisms}, {\it isomorphisms}, {\it subalgebras}
of open-string vertex algebras are defined in the obvious way.}
\end{defn}

We shall denote the open-string vertex algebra by $(V, Y^{O}, 
\one, D)$ or simply $V$. For $u\in V$ and $r\in \R_{+}$, 
we call the map $Y^{O}(u, r): V\to \overline{V}$ the vertex 
operator associated to $u$ and $r$. 

\begin{rema}\label{perm-p}
{\rm Note that in the definition above, the real number $r$ in the vertex 
operator $Y^{O}(u, r)$ is positive, not in $\R^{\times}$. So a natural 
question is whether one has natural vertex operators associated to 
negative real numbers so that we have a vertex map $Y^{O}$ {}from
$(V\otimes V) \times \R^{\times}$ to $\overline{V}$. The answer is yes. 
For any $u, v\in V$ and $r\in 
-\R_{+}$, we define
\begin{equation}\label{permute}
Y^{O}(u, r)v=e^{rD}Y^{O}(v, -r)u.
\end{equation}
(Note that $e^{rD}Y^{O}(v, -r)u$ is a well-defined element of 
$\overline{V}$ by the local-truncation property for $D'$.)
Note that (\ref{permute}) resembles
the skew-symmetry for vertex operator algebras. We know that 
the skew-symmetry is analogous to  commutativity for commutative
associative algebras.  But (\ref{permute})
does not give a skew-symmetry property and is not an 
analogue of the commutativity mentioned above. 
Instead, (\ref{permute}) is an analogue of the 
relation between the product and the opposite product for 
an associative algebra. In fact, for
an associative algebra $V$, we can define an opposite product
\begin{equation}\label{prod-op-prod}
(uv)^{{\rm op}}=vu
\end{equation}
for $u, v\in V$.
We can also define an open-string
vertex algebra in terms of a vertex map of the form 
$V\otimes V\times \R^{\times}\to \overline{V}$ and then 
(\ref{permute}) becomes an axiom. In applications, it is convenient to have
vertex operators associated 
to negative numbers. For example, since
\begin{eqnarray*}
\langle v', Y^{O}(Y^{O}(v_{1}, r_0)v_{2}, r_2)v\rangle
&=&\sum_{m\in \R}\langle v', Y^{O}(P_{m}Y^{O}(v_{1}, r_0)v_{2}, r_2)v\rangle\nn
&=&\sum_{m\in \R}\langle e^{r_{2}D'}v', Y^{O}(v, -r_2)P_{m}Y^{O}(v_{1}, r_0)
v_{2}\rangle,
\end{eqnarray*}
the left-hand side (a matrix elements of an iterate of vertex operators) is 
absolutely convergent when $r_{2}>r_{0}>0$ 
if and only if the right-hand side (a matrix element of a product of 
vertex operators) is. (Note that 
by the local-truncation property for $D'$,
$e^{r_{2}D'}v'\in V'$.)
In fact, one can prove that 
for $v_{1}, \dots, v_{n}, v\in V$
and $v'\in V'$, the series 
\begin{eqnarray*}
\lefteqn{\langle v', Y^{O}(v_{1}, r_1)\cdots Y^{O}(v_{n}, r_n)v\rangle}\nn
&&=\sum_{m_{1}, \dots, m_{n-1}\in \R}
\langle v', Y^{O}(v_{1}, r_1)P_{m_{1}}Y^{O}(v_{2}, r_2)\cdots 
P_{m_{n-1}}Y^{O}(v_{n}, r_n)v\rangle
\end{eqnarray*}
converges absolutely 
when $|r_1|>\cdots >|r_n|>0$. One can also prove 
the absolute convergence of all the products and iterates of 
vertex operators associated to real numbers in natural regions. 
For the skew-symmetry
for $Y^{O}$, see Remark \ref{skew-sym}.}
\end{rema}

We still use $\mathbf{d}$ to denote the natural 
extension of $\mathbf{d}$ to an element of $\hom(\overline{V},
\overline{V})$.

\begin{prop}  
The $\mathbf{d}$-bracket property (\ref{d-bra}) for
all $u\in V$ and $r\in \R_{+}$,
is equivalent to the {\it $\mathbf{d}$-conjugation 
property}
\begin{equation} \label{d-conj}
a^{\mathbf{d}}Y^{O}(u,r)a^{-\mathbf{d}}=Y^{O}(a^{\mathbf{d}}u,ar)
\end{equation}
for all $u\in V$, $r\in \R_{+}$ and $a\in \R_{+}$.
We also have $\one\in V_{(0)}$ and 
$D\one=0$. 
\end{prop}
\pf
If (\ref{d-conj}) holds for all 
$u\in V$, $r\in \R_{+}$ and $a\in \R_{+}$,
then 
\begin{equation}\label{d-conj-1}
e^{s\mathbf{d}}Y^{O}(u,r)e^{-s\mathbf{d}}=Y^{O}(e^{s\mathbf{d}}u,e^{s}r)
\end{equation}
for all $u\in V$, $r\in \R_{+}$ and $s\in \R$.
Taking the derivative with respect to $s$ of 
both sides of (\ref{d-conj-1})
and then letting $s=0$, we obtain (\ref{d-bra}).

Conversely, assume that (\ref{d-bra}) holds for all 
$u\in V$ and $r\in \R_{+}$. 
Let $u,v\in V$ and $v'\in V'$ be homogeneous. We have
\begin{eqnarray}\label{d-conj-2}
\langle v', [\mathbf{d}, Y^{O}(u,r)]v\rangle &=&
\langle \mathbf{d}'v', , Y^{O}(u,r)v\rangle 
-\langle v', , Y^{O}(u,r)\mathbf{d} v\rangle\nn
&=&(\wt v'-\wt v)
\langle w, Y^{O}(u,r)v\rangle,
\end{eqnarray}
where $\mathbf{d}'$ is the adjoint of $\mathbf{d}$.
On the other hand, 
\begin{equation}\label{d-conj-3}
\left\langle v', \left( Y^{O}(\mathbf{d} u,r) +
r\frac{d}{dr} Y^{O}(u,r) \right) v\right\rangle
= \left(\wt u + r\frac{d}{dr}\right) \langle v', Y^{O}(u,r)v\rangle.
\end{equation}
By (\ref{d-bra}), (\ref{d-conj-2}) and (\ref{d-conj-3}), 
we see that $f(r)=\langle v', Y^{O}(u,r)v\rangle$ satisfies the 
differential equation
$$
r\frac{df(r)}{dr} = (\wt v'-\wt u-\wt v) f(r)
$$
Any solution of this equation 
is of the form $Cr^{\swt w-\swt u-\swt v}$
for some $C\in \C$. In particular, $\langle v', Y^{O}(u,r)v\rangle$
is of this form. Therefore,
\begin{eqnarray*}
\langle v', a^{\mathbf{d}}Y^{O}(u,r)a^{-\mathbf{d}} v \rangle
&=&\langle a^{\mathbf{d}}v', Y^{O}(u,r)a^{-\mathbf{d}} v \rangle\nn
&=&a^{\swt v'-\swt v}\langle v', Y^{O}(u,r)v\rangle\nn
&=&Ca^{\swt v'-\swt v}r^{\swt w-\swt u-\swt v}\nn
&=&Ca^{\swt u}(ar)^{\swt w-\swt u-\swt v}\nn
&=& \langle v', Y^{O}(a^{\mathbf{d}}u, ar)v\rangle.
\end{eqnarray*}
Since such $u, v$ span $V$ and such $v'$ spans 
$V'$, we obtain (\ref{d-conj}).

The identity property, 
the creation property and (\ref{d-bra})
imply $\mathbf{d}\one=0$ which means $\one\in V_{(0)}$. 
The identity property and the 
$D$-derivative property imply $D\one=0$. 
\epfv

The $\mathbf{d}$-conjugation property also has the following 
very important consequence:

\begin{prop}\label{expn}
For $u\in V$, there exist $u^{+}_{n}\in 
\edo V$ of weights $\wt u-n-1$ for $n\in \R$ such that for $r\in \R_{+}$,
\begin{equation}\label{expn-1}
Y^{O}(u, r)=\sum_{n\in \R}u^{+}_{n}r^{-n-1}.
\end{equation}
\end{prop}
\pf
For homogeneous $u\in V$ and $n\in \R$, 
let $u^{+}_{n}\in \edo V$ be defined by
$$u^{+}_{n}v=P_{\swt u-n-1+\swt v}Y^{O}(u, 1)v$$
for homogeneous $v\in V$. Then by the $\mathbf{d}$-conjugation property,
for any homogeneous $u, v\in V$,
\begin{eqnarray*}
Y^{O}(u, r)v&=&r^{\mathbf{d}}Y^{O}(r^{-\mathbf{d}}u, 1)r^{-\mathbf{d}}v\nn
&=&r^{-\swt u-\swt v}r^{\mathbf{d}}\sum_{n\in \R}
P_{\swt u-n-1+\swt v}Y^{O}(u, 1)v\nn
&=&\sum_{n\in \R}P_{\swt u-n-1+\swt v}Y^{O}(u, 1)vr^{-n-1}\nn
&=&\sum_{n\in \R}u^{+}_{n}vr^{-n-1}. 
\end{eqnarray*}
\epf

\begin{rema}
{\rm In the proposition above, (\ref{expn-1}) holds only for 
$r\in \R_{+}$. In fact, there are also $u^{-}_{n}\in 
\edo V$ of weights $\wt u-n-1$ for $n\in \R$ such that 
for $r\in -\R_{+}$,
$$Y^{O}(u, r)=\sum_{n\in \R}u^{-}_{n}r^{-n-1}.$$
But in general $u^{-}_{n}\ne u^{+}_{n}$. In this paper, we shall not use
$u_{n}^{-}$, $n\in \R$.}
\end{rema}

{}From Proposition \ref{expn}, we see that for any $u\in V$, there is a 
formal-variable vertex operator 
$$\Y^{f}(u, x)=\sum_{n\in \R}u^{+}_{n}x^{-n-1}\in (\edo V)[[x, x^{-1}]]$$
where $x$ is a formal variable. 
We shall also use the notation 
$\Y^{f}(u, z)$ to denote the vertex operator
associated to $u\in V$ and a nonzero complex number $z$, that
is,
$$\Y^{f}(u, z)=\sum_{n\in \R}u^{+}_{n}z^{-n-1}.$$
(Note that by our convention, for $z\in \C^{\times}$, 
$z^{-n-1}=e^{(-n-1)\log z}$ for $n\in 
\R$, where $\log z=\log |z|+i\arg z$, $0\le \arg z<2\pi$.) 
Thus for  $u\in V$, $Y^{O}(u, r)=\Y^{f}(u, r)$ for $r\in \R_{+}$ but 
in general $Y^{O}(u, r)\ne \Y^{f}(u, r)$ for $r\in -\R_{+}$.

\begin{rema}\label{skew-sym}
{\rm As we have discussed in Remark \ref{perm-p}, 
(\ref{permute}) has nothing to do with skew-symmetry. In fact, if 
$\Y^{f}$ satisfies 
$$\Y^{f}(u, x)v=e^{xL(-1)}\Y^{f}(v, y)u\lbar_{y^{n}=e^{n\pi i}x^{n},\;
n\in \R}$$
for $u, v\in V$, then we say that $Y^{O}$ has {\it skew-symmetry}. 
(For simplicity, in the remaining part of this paper, 
we shall use $\Y^{f}(v, -x)u$ to denote 
$\Y^{f}(v, y)u|_{y^{n}=e^{n\pi i}x^{n},\;
n\in \R}$.)
Note that skew-symmetry for $Y^{O}$ gives a relation between 
$Y^{O}(u, r)v$ 
and its analytic extension to the negative real line 
for $u, v\in V$ and $r\in \R_{+}$
while (\ref{permute}) gives a relation between 
$Y^{O}(u, r)v$ and $Y^{O}(v, -r)u$ for $u, v\in V$ and $r\in \R_{+}$.
Clearly, these two relations are in general different.}
\end{rema}

\begin{prop}\label{formal-conj}
The $\mathbf{d}$-bracket  and $D$-derivative properties
hold for $\Y^{f}$, that is, (\ref{d-bra}) and (\ref{D-der})
hold when $Y^{O}$ is replaced by $\Y^{f}$ and 
$r$ is replaced by the formal variable $x$. 
We also have the following {$\mathbf{d}$-} and {$D$-conjugation properties}:
For $u\in V$ and $y$ another formal variable,
$$y^{\mathbf{d}}\Y^{f}(u, x)y^{-\mathbf{d}}=\Y^{f}(y^{\mathbf{d}}u,yx),$$
and 
\begin{equation} \label{D-conj}
\Y^{f}(u, x+y)=e^{yD}\Y^{f}(u, x)e^{-yD}=\Y^{f}(e^{yD}u, x).
\end{equation}
In particular, these conjugation formulas also hold when we 
substitute suitable complex numbers for $x$ and $y$ such that 
both sides of these formulas make sense  as (or converges to) 
maps {}from $V$ to $\overline{V}$.
\end{prop}
\pf
The $\mathbf{d}$-bracket formula and $D$-derivative property 
follow {}from the definition of the formal-variable vertex operators
and the corresponding properties for the defining vertex map. 
The $\mathbf{d}$-conjugation property for the formal vertex operator 
follows immediately {}from (\ref{d-conj}). 
The $D$-conjugation property follows {}from the 
$D$-derivative property.
\epfv

We have the following easy consequence of Proposition \ref{formal-conj}:

\begin{cor}
For any $u\in V$, 
\begin{equation}\label{creation}
\Y^{f}(u, x)\one=e^{xD}u.
\end{equation}
\end{cor}
\pf
By the creation property, we see that for any $r\in \R_{+}$ and any $u\in V$, 
$$Y^{O}(u, r)\one=\sum_{n\in (-\R_{+}-1)\cup \{-1\}}u_{n}^{+}\one r^{-n-1}$$
and $u_{-1}^{+}\one=u$. But by the $D$-derivative property,
\begin{eqnarray*}
\lim_{r\to 0}\frac{d^{k}}{dr^{k}}Y^{O}(u, r)\one
&=&\lim_{r\to 0}Y^{O}(D^{k}u, r)\one\nn
&=&D^{k}u
\end{eqnarray*}
for $k\in \N$. Thus we see that 
$$Y^{O}(u, r)\one=\sum_{n\in -\Z_{+}}u_{n}^{+}\one r^{-n-1}.$$
So $\Y^{f}(u, x)\one$
is a power series in $x$ and $\lim_{x\to 0} \Y^{f}(u, x)\one=u$. 
By these properties, $D\one=0$ and 
the $D$-conjugation property for $\Y^{f}$, we obtain
\begin{eqnarray*}
\Y^{f}(u, y)\one&=&\lim_{x\to 0}\Y^{f}(u, x+y)\one\nn
&=&\lim_{x\to 0}e^{yD}\Y^{f}(u, x)e^{-yD}\one\nn
&=&\lim_{x\to 0}e^{yD}\Y^{f}(u, x)\one\nn
&=&e^{yD}u,
\end{eqnarray*}
proving (\ref{creation}).
\epf

\begin{prop}\label{assoc-y-f}
The formal vertex operator map $\Y^{f}$ has the following 
properties:

\begin{enumerate}

\item Convergence: The series 
\begin{eqnarray}
&\langle v', \Y^{f}(v_{1}, z_{1})\Y^{f}(v_{2}, z_{2})v\rangle,&\label{prod1}\\
&\langle v', \Y^{f}(v_{2}, z_{2})\Y^{f}(v_{1}, z_{1})v\rangle,&\label{prod2}\\
&\langle v', \Y^{f}(\Y^{f}(v_{1}, z_{1}-z_{2})v_{2}, z_{2})v\rangle,&
\label{iter1}\\
&\langle v', \Y^{f}(\Y^{f}(v_{2}, z_{2}-z_{1})v_{1}, z_{1})v\rangle&
\label{iter2}
\end{eqnarray}
are absolutely convergent in the regions $|z_{1}|>|z_{2}|>0$, 
$|z_{2}|>|z_{1}|>0$, $|z_{2}|>|z_{1}-z_{2}|>0$, $|z_{1}|>|z_{1}-z_{2}|>0$,
respectively.

\item Associativity: For 
$v_{1}, v_{2}, v\in V$ and $v'\in V'$,
(\ref{prod1}) and (\ref{iter1}) are equal in the region 
$|z_{1}|>|z_{2}|>|z_{1}-z_{2}|>0$, and (\ref{prod2}) and (\ref{iter2})
are equal in the region $|z_{2}|>|z_{1}|>|z_{1}-z_{2}|>0$. 

\end{enumerate}

\end{prop}
\pf
By definition, 
(\ref{prod1}), (\ref{prod2}), (\ref{iter1}) and (\ref{iter2})
converge absolutely when $z_{1}, z_{2}\in \R_{+}$ satisfying 
$z_{1}>z_{2}>0$, $z_{2}>z_{1}>0$, $z_{2}>z_{1}-z_{2}>0$ and
$z_{1}>z_{1}-z_{2}>0$,
respectively. Consequently, (\ref{prod1}),(\ref{prod2}),
(\ref{iter1}) and (\ref{iter2}) converge absolutely 
for $z_{1}, z_{2}\in \C$ satisfying $|z_{1}|>|z_{2}|>0$, 
$|z_{2}|>|z_{1}|>0$, $|z_{2}|>|z_{1}-z_{2}|>0$ $|z_{1}|>|z_{1}-z_{2}|>0$,
respectively. The convergence is
proved.

In particular, (\ref{prod1})
and (\ref{iter1}) give (possibly multivalued) 
analytic functions defined 
on the regions $|z_{1}|>|z_{2}|>0$ and
$|z_{2}|>|z_{1}-z_{2}|>0$, respectively.
By associativity for $Y^{O}$, (\ref{prod1}) and (\ref{iter1})
are equal for $z_{1}, z_{2}\in \R_{+}$ satisfying 
$z_{1}>z_{2}>z_{1}-z_{2}>0$. By the basic properties of analytic
functions, (\ref{prod1}) and (\ref{iter1})
are equal for $z_{1}, z_{2}\in \C$ satisfying 
$|z_{1}|>|z_{2}|>|z_{1}-z_{2}|>0$ (the intersection of the 
regions $|z_{1}|>|z_{2}|>0$ and $|z_{2}|>|z_{1}-z_{2}|>0$
on which the analytic functions 
(\ref{prod1}) and (\ref{iter1}) are defined). The 
second part of the associativity for $\Y^{f}$
can be obtained {}from the 
first part by substituting  
$v_{2}$, $v_{1}$, $z_{2}$ and $z_{1}$ for
$v_{1}$, $v_{2}$, $z_{1}$ and $z_{2}$.
\epf

\begin{defn} 
{\rm A {\it grading-restricted open-string vertex 
algebra} is an open-string vertex algebra satisfying the 
following conditions:

\bnu

\setcounter{enumi}{7}

\item The {\it grading-restriction conditions}: For all $n\in \R$, 
$\dim V_{(n)}<\infty$ (the {\it finite-dimensionality of 
homogeneous subspaces}) and $V_{(n)}=0$ when 
$n$ is sufficiently negative (the {\it lower-truncation 
condition for grading}).

\enu

A {\it  conformal open-string  vertex algebra} 
is an open-string vertex algebra equipped with a {\it conformal element}
$\omega\in V$ satisfying the following 
conditions:

\bnu

\setcounter{enumi}{8}

\item The {\it Virasoro relations}: For any $m, n\in \Z$, 
$$[ L(m), L(n)]  = (m-n)L(m+n) - \frac{c}{12}(m^3-m)\delta_{m+n,0},$$
where $L(n)$, $n\in \Z$ are given by
$$Y^{O}(\omega, r) =\sum_{n\in \Z} L(n) r^{-n-2}$$
and $c\in \C$.

\item The {\it commutator formula for Virasoro operators and 
formal vertex operators (or component operators)}: For $v\in V$,
$\Y^{f}(\omega, x)v$ involves only finitely many negative powers
of $x$ and
$$[\Y^{f}(\omega, x_{1}), \Y^{f}(v, x_{2})]
=\res_{0}x_{2}^{-1}\delta\left(\frac{x_{1}-x_{0}}{x_{2}}\right)
\Y^{f}(\Y^{f}(\omega, x_{0})v, x_{2}).$$

\item The {\it $L(0)$-grading property} and 
{\it $L(-1)$-derivative property}:  $L(0)=\mathbf{d}$ and $L(-1)=D$.

\enu   

A {\it grading-restricted conformal open-string vertex algebra}
or {\it open-string vertex operator algebra}
is a conformal open-string vertex algebra satisfying the 
grading-restriction condition.}
\end{defn}

We shall denote the conformal open-string vertex algebra defined above
by $(V, Y^{O}, \one, \omega)$ or simply $V$. The complex number $c$ 
in the definition is called the {\it central charge} 
of the algebra. Note that 
the grading-restriction conditions imply the local-truncation
property for $D'$.

\begin{prop}  
Let $V$ be a grading-restricted open-string vertex algebra.
Then for $u,v\in V$, $u^{+}_nv=0$ if $n$ is sufficiently negative. 
\end{prop}
\pf
This follows immediately {}from 
the lower-truncation 
condition for grading and the fact that the weights of 
$u^{+}_{n}$ for $n\in \R$ is $\wt u-n-1$.
\epf

\renewcommand{\theequation}{\thesection.\arabic{equation}}
\renewcommand{\thethm}{\thesection.\arabic{thm}}
\setcounter{equation}{0}
\setcounter{thm}{0}

\section{Intertwining operators and open-string vertex algebras}

In this section, we establish a connection between open-string 
vertex algebras and  intertwining operator algebras. 
We assume that the reader is familiar with the basic 
notions and properties in the representation theory
of vertex operator algebras and we also assume that the reader is
familiar with the notion of intertwining 
operator algebra. See \cite{FHL}, \cite{H6} and \cite{H7} for details.

Let $V$ be an open-string vertex algebra and 
$S$ a subset of $V$. Then the 
{\it open-string vertex subalgebra of $V$ generated by $S$}
is the smallest open-string vertex subalgebra of $V$ containing $S$.

\begin{prop}
Let $V$ be a conformal open-string vertex algebra and 
$\langle \omega \rangle$ the open-string vertex subalgebra of $V$ generated
by $\omega$.   Then  $\langle \omega \rangle$ is in fact a vertex operator
algebra. In particular, $V$ is a module for the 
vertex operator algebra $\langle \omega \rangle$.
\end{prop}
\pf
All the axioms for a vertex operator algebra are satisfied by 
$\langle \omega \rangle$ obviously except for the commutativity
or equivalently the commutator formula. But the Virasoro relations
imply the commutator formula for the vertex operators 
for $\langle \omega \rangle$.
\epfv

More generally, we have the following generalization:
Let $V$ be an open-string vertex algebra and let 
\begin{eqnarray*}
\lefteqn{C_{0}(V)=\Bigg\{ u\in \coprod_{n\in \Z}V_{(n)}
\;\Big|\;\Y^{f}(u, x)\in (\edo V)[[x, x^{-1}]],}\nn
&&\quad\quad\quad\quad\quad\quad\quad\quad\quad\quad\quad\quad
\Y^{f}(v, x)u=e^{xD}\Y^{f}(u, -x)v,\;\forall v\in V \Bigg\}.
\end{eqnarray*}
In particular, for elements of $C_{0}(V)$, skew-symmetry holds.
Clearly $C_{0}(V)$ is not zero since by (\ref{creation}),
$\one \in C_0(V)$.

For an open-string vertex algebra $V$,
the formal vertex operator map $\Y^{f}$ for $V$ induces a map {}from
$C_{0}(V)\otimes C_{0}(V)$ to $V[[x, x^{-1}]]$. We denote this map by 
$\Y^{f}|_{C_{0}(V)}$. We first need:

\begin{prop}\label{anal-ext}
Let $v_{1}\in C_{0}(V)$, $v_{2}, v\in V$ and $v'\in V'$. Then 
there exists a (possibly multivalued) analytic function on 
$$M^{2}=\{(z_{1}, z_{2})\in \C^{2}\;|\; z_{1}, 
z_{2}\ne 0, z_{1}\ne z_{2}\}$$
such that it  is single valued in $z_{1}$ and 
is equal to the (possibly multivalued) analytic extensions of 
(\ref{prod1}), (\ref{prod2}), (\ref{iter1}) and (\ref{iter2})
in the regions $|z_{1}|>|z_{2}|>0$, 
$|z_{2}|>|z_{1}|>0$, $|z_{2}|>|z_{1}-z_{2}|>0$ and 
$|z_{1}|>|z_{1}-z_{2}|>0$, respectively. Moreover, if 
$v_{2}$ is in $C_{0}(V)$, then this analytic function 
is single valued in both $z_{1}$ and $z_{2}$. If 
$V$ satisfies the grading-restriction condition,
then this analytic function is a rational function 
with the only possible poles $z_{1}, z_{2}=0$ and $z_{1}=z_{2}$.
\end{prop}
\pf
By Proposition \ref{assoc-y-f}, 
(\ref{prod1}), (\ref{prod2}) and (\ref{iter1})
are absolutely convergent
in the regions $|z_{1}|>|z_{2}|>0$, 
$|z_{2}|>|z_{1}|>0$, $|z_{2}|>|z_{1}-z_{2}|>0$, respectively, and 
the associativity for $\Y^{f}$ holds. 

Since $v_{1}\in C_{0}(V)$, by definition,
$\Y^{f}(v_{1}, x)v_{2}\in V[[x, x^{-1}]]$ and 
we have the skew-symmetry
\begin{eqnarray*}
\Y^{f}(v_{1}, x)v_{2}&=&e^{xD}\Y^{f}(v_{2}, -x)v_{1},\\
\Y^{f}(v_{2}, x)v_{1}&=&e^{xD}\Y^{f}(v_{1}, -x)v_{2}.
\end{eqnarray*}
In \cite{H6} it was proved that commutativity for intertwining operators
follows {}from associativity and skew-symmetry for 
intertwining operators. For reader's 
convenience, here we give a proof of commutativity 
in the special case
in which we are interested.  

By associativity, (\ref{prod1}) and (\ref{iter1}) are equal in the region
$|z_{1}|>|z_{2}|>|z_{1}-z_{2}|>0$. By associativity also, 
(\ref{prod2}) and (\ref{iter2})
converge absolutely to analytic functions
defined on the regions $|z_{2}|>|z_{1}|>0$ and $|z_{1}|>|z_{1}-z_{2}|>0$,
respectively, and are equal in the region
$|z_{2}|>|z_{1}|>|z_{1}-z_{2}|>0$.
By skew-symmetry and the $D$-derivative 
property, for $z_{1}, z_{2}\in \mathbb{C}$ satisfying 
$|z_{1}|>|z_{1}-z_{2}|>0$ and $|z_{2}|>|z_{1}-z_{2}|>0$,
we have 
\begin{eqnarray*}
\lefteqn{\langle v', \Y^{f}(\Y^{f}(v_{1}, z_1-z_2)v_{2},z_2)v\rangle}\nn
&&=
\langle v', \Y^{f}(e^{(z_{1}-z_{2})D}
\Y^{f}(v_{2}, -(z_1-z_2))v_{1}, z_2)v\rangle\nn
&&=\langle v', \Y^{f}(\Y^{f}(v_{2}, z_2-z_1)v_{1}, z_{2}+(z_{1}-z_2))v
\rangle\nn
&&=\langle v', \Y^{f}(\Y^{f}(v_{2}, z_2-z_1)v_{1}, z_{1})v\rangle,
\end{eqnarray*}
that is, in the region given by
$|z_{1}|>|z_{1}-z_{2}|>0$ and $|z_{2}|>|z_{1}-z_{2}|>0$,
(\ref{iter1}) and (\ref{iter2}) are equal. 
Since 
(\ref{prod1}) is equal to (\ref{iter1}) in the region
$|z_{1}|>|z_{2}|>|z_{1}-z_{2}|>0$, (\ref{iter1}) is equal to 
(\ref{iter2}) in the region given by
$|z_{1}|>|z_{1}-z_{2}|>0$ and $|z_{2}|>|z_{1}-z_{2}|>0$, and 
(\ref{iter2}) is equal to (\ref{prod2}) in the region 
$|z_{2}|>|z_{1}|>|z_{1}-z_{2}|>0$, we see that 
(\ref{prod1}) and (\ref{prod2}) are analytic extensions of 
each other. So commutativity is proved. 

Now we prove the existence of the function stated in the proposition. 
By skew-symmetry, we have 
$$\Y^{f}(v, z)\one=e^{zD}\Y^{f}(\one, -z)v=e^{zD}v$$
for any $v\in C_{0}(V)$. Thus by definition, 
for $v_{1}\in C_{0}(V)$, $v_{2}, v\in V$
and $v'\in (C_{0}(V))'$,
$$\langle v', \Y^{f}(v_{1}, z_1)\Y^{f}(v_{2}, z_2)e^{z_{3}D}v\rangle
=\langle v', \Y^{f}(v_{1}, z_1)\Y^{f}(v_{2}, z_2)\Y^{f}(v, z_{3})\one\rangle$$
converges absolutely for $z_{1}, z_{2}, z_{3}\in \R^{\times}$ satisfying
$|z_{1}|>|z_{2}|>|z_{3}|>0$. Consequently it also converges absolutely
for $z_{1}, z_{2}, z_{3}\in \C$ satisfying
$|z_{1}|>|z_{2}|>|z_{3}|>0$. Now the same proof as the one for 
Lemma 4.1 in \cite{H6} shows that there exists a 
(possibly multivalued) analytic functions on $M^{2}$ 
such that it is equal to (possibly multivalued) analytic extensions of
(\ref{prod1}), (\ref{prod2}), (\ref{iter1}) and (\ref{iter2})
in the regions $|z_{1}|>|z_{2}|>0$, 
$|z_{2}|>|z_{1}|>0$, $|z_{2}|>|z_{1}-z_{2}|>0$ and 
$|z_{1}|>|z_{1}-z_{2}|>0$, respectively.
Since (\ref{prod1}), (\ref{prod2}) and (\ref{iter1})
give analytic functions which are all single valued in $z_{1}$, 
this function as the analytic extension of these functions
must also be 
single valued in $z_{1}$.

If $v_{2}$ is in $C_{0}(V)$, then by definition,
$\Y^{f}(v_{2}, x)v\in V[[x, x^{-1}]]$ and thus 
(\ref{prod1}), (\ref{prod2}) and (\ref{iter2}) 
give analytic functions which are also single valued in $z_{2}$.
So their analytic extension
is also single valued in both $z_{1}$ and $z_{2}$. 
If $V$ 
satisfies the grading-restriction condition, then the 
singularities $z_{1}, z_{2}=0, \infty$ and $z_{1}=z_{2}$
of this analytic extension 
are all poles and is therefore 
a rational function in $z_{1}$ and $z_{2}$ with the 
only possible poles $z_{1}, z_{2}=0$ and $z_{1}=z_{2}$.
\epf

\begin{thm}\label{op-va-int} 
Let $V$ be a  grading-restricted open-string vertex algebra. Then 
the image of $C_{0}(V)\otimes C_{0}(V)$ under
$\Y^{f}|_{C_{0}(V)}$ is in $C_{0}(V)[[x, x^{-1}]]$
and the image of $C_{0}(V)$ under $D$ is in $C_{0}(V)$. Moreover,
$$(C_{0}(V), \Y^{f}|_{C_{0}(V)}, \one, D)$$
is a grading-restricted vertex  algebra, 
$V$ is a $C_0(V)$-module and 
$\Y^{f}$ is an intertwining operator of type $\binom{V}{VV}$ 
for the vertex algebra $C_{0}(V)$. 
\end{thm}
\pf
Let $v_{1}, v_{2}$ be homogeneous elements of $C_{0}(V)$. 
We would like to show that 
$\Y^{f}(v_{1}, x)v_{2}\in C_{0}(V)[[x, x^{-1}]]$. 
First of all, since $v_{1}\in C_{0}(V)$,
$\Y^{f}(v_{1}, x)v_{2}\in V[[x, x^{-1}]]$. Since 
$v_{1}, v_{2}\in C_{0}(V)$, $\wt v_{1}, \wt v_{2}\in \Z$. 
Thus by Proposition \ref{expn}, 
$\Y^{f}(v_{1}, x)v_{2}\in 
\left(\coprod_{n\in \Z}V_{(n)}\right)[[x, x^{-1}]]$.
By Proposition \ref{anal-ext}, the analytic extension of
(\ref{iter1}) to $M^{2}$ is a single-valued analytic function.
In particular, (\ref{iter1}) gives a single-valued  analytic 
function in $z_{1}$ and $z_{2}$. 
Thus 
$$\Y^{f}(\Y^{f}(v_{1}, x)v_{2}, x_{2})v
\in (V[[x_{2}, x_{2}^{-1}]])[[x, x^{-1}]].$$

For $v\in V$, $v'\in V'$ and $z_{1}, z_{2}\in \R_{+}$ satisfying 
$z_{1}>z_{2}>z_{1}-z_{2}>0$, 
\begin{eqnarray}\label{mero-c-1}
\langle v', \Y^{f}(\Y^{f}(v_{1}, z_{1}-z_{2})v_{2}, z_{2})v\rangle
&=&\langle v', \Y^{f}(v_{1}, z_{1})\Y^{f}(v_{2}, z_{2})v\rangle\nn
&=&\langle v', \Y^{f}(v_{1}, z_{1})e^{z_{2}D}\Y^{f}(v, -z_{2})v_{2}
\rangle.\nn
&&
\end{eqnarray}
The right-hand side of (\ref{mero-c-1}) is well defined 
when $z_{1}, z_{2}\in \C$ and
$|z_{1}|>|z_{2}|>0$ and is equal to 
\begin{eqnarray}\label{mero-c-2}
\lefteqn{\langle v', e^{z_{2}D}\Y^{f}(v_{1}, z_{1}-z_{2})
\Y^{f}(v, -z_{2})v_{2}\rangle}\nn
&&=\langle v', e^{z_{1}D}\Y^{f}(\Y^{f}(v, -z_{2})v_{2}, -(z_{1}-z_{2}))v_{1}
\rangle\nn
&&=\langle v', e^{z_{1}D}\Y^{f}(v, -z_{1})\Y^{f}(v_{2}, -(z_{1}-z_{2}))v_{1}
\rangle\nn
&&=\langle v',  e^{z_{1}D}\Y^{f}(v, -z_{1})e^{-(z_{1}-z_{2})D}
\Y^{f}(v_{1}, z_{1}-z_{2})v_{2}\rangle
\end{eqnarray}
when $z_{1}, z_{2}\in \R_{+}$ and $z_{1}>z_{1}-z_{2}>z_{2}>0$.
The right-hand side of (\ref{mero-c-2}) is well defined when 
$z_{1}, z_{2}\in \C$ and $|z_{1}|>|z_{1}-z_{2}|>0$ and is equal to 
\begin{equation}\label{mero-c-3}
\langle v',  e^{z_{2}D}\Y^{f}(v, -z_{2})
\Y^{f}(v_{1}, z_{1}-z_{2})v_{2}\rangle
\end{equation}
when $z_{1}, z_{2}\in \C$ and $|z_{1}|>|z_{2}|>|z_{1}-z_{2}|>0$.
{}From (\ref{mero-c-1})--(\ref{mero-c-3}), we see that 
the left-hand side of (\ref{mero-c-1}) and the right-hand side of
(\ref{mero-c-3}) are analytic extensions of each other. 
Since both the left-hand side of (\ref{mero-c-1}) and the right-hand side 
of (\ref{mero-c-3}) are well defined single-valued analytic 
functions on the region $|z_{2}|>|z_{1}-z_{2}|>0$,
they are equal when $|z_{2}|>|z_{1}-z_{2}|>0$. Thus 
we obtain 
$$\Y^{f}(\Y^{f}(v_{1}, x)v_{2}, x_{2})v
=e^{x_{2}D}\Y^{f}(v, -x_{2})
\Y^{f}(v_{1}, x)v_{2}$$
where $x$ and $x_{2}$ are two commuting formal variables. 
So $\Y^{f}(v_{1}, x)v_{2}\in C_{0}(V)[[x, x^{-1}]]$.

Let $u$ be a homogeneous element of $C_{0}(V)$. Then 
$\wt u\in \Z$. Since $D$ has weight $1$, $Du\in \coprod_{n\in \Z}V_{(n)}$.
By the $D$-derivative property, we see that 
$$\Y^{f}(Du, x)=\frac{d}{dx}\Y^{f}(u, x)\in (\edo V)[[x, x^{-1}]].$$
For any $v\in V$, using the $D$-derivative property and the 
$D$-bracket formula, we obtain
\begin{eqnarray*}
\Y^{f}(Du, x)v&=&\frac{d}{dx}\Y^{f}(u, x)v\nn
&=&\frac{d}{dx}e^{xD}\Y^{f}(v, -x)u\nn
&=&e^{xD}D\Y^{f}(v, -x)u-e^{xD}\Y^{f}(Dv, -x)u\nn
&=&e^{xD}\Y^{f}(v, -x)Du.
\end{eqnarray*}
So $Du\in C_{0}(V)$.

To show that $C_{0}(V)$ is a vertex algebra,
we need only verify commutativity, associativity and 
rationality since all the other axioms are clearly satisfied.
But associativity, commutativity and rationality 
have been proved in Proposition 
\ref{anal-ext}. The proof of the fact that $V$ is a $C_{0}(V)$-module and 
$\Y^{f}$ is an intertwining operator of type $\binom{V}{VV}$ 
for $C_{0}(V)$ is completely the 
same. 
\epf

\begin{prop}
Let $V$ be a conformal open-string
vertex algebra. Then $\omega\in C_0(V)$. 
\end{prop}
\pf
By definition, $\omega\in \coprod_{n\in \Z}V_{(n)}$
and $\Y^{f}(\omega, x)\in (\edo V)[[x, x^{-1}]]$. 
For any $v\in V$, the commutator formula for $\omega$
and formal vertex operators implies the commutativity
for $\Y^{f}(\omega, z_{1})$ and $\Y^{f}(v, z_{2})$. 
In particular, for any $v'\in V'$, 
\begin{equation}\label{prod1-omega}
\langle v', \Y^{f}(\omega, z_{1})\Y^{f}(v, z_{2})\one\rangle
\end{equation}
and 
\begin{equation}\label{prod2-omega}
\langle v', \Y^{f}(v, z_{2})\Y^{f}(\omega, z_{1})\one\rangle
\end{equation}
are absolutely convergent in the regions $|z_{1}|>|z_{2}|>0$ and 
$|z_{2}|>|z_{1}|>0$, respectively, and are analytic extensions
of each other. Also by associativity we know that 
\begin{equation}\label{iter-omega}
\langle v', \Y^{f}(\Y^{f}(\omega, z_{1}-z_{2})v, z_{2})\one\rangle
\end{equation}
and 
\begin{equation}\label{iter2-omega}
\langle v', \Y^{f}(\Y^{f}(v, z_{2}-z_{1})\omega, z_{1})\one\rangle
\end{equation}
are absolutely convergent in the region $|z_{2}|>|z_{1}-z_{2}|>0$ 
and $|z_{1}|>|z_{1}-z_{2}|>0$, respectively, 
and are equal to (\ref{prod1-omega}) and (\ref{prod2-omega}),
respectively, in the region 
$|z_{1}|>|z_{2}|>|z_{1}-z_{2}|>0$ and 
$|z_{2}|>|z_{1}|>|z_{1}-z_{2}|>0$, respectively. Thus 
(\ref{iter-omega}) and
(\ref{iter2-omega}) are also analytic extensions of each other.
Note that by (\ref{D-conj}),
\begin{eqnarray}\label{iter3-omega}
\lefteqn{\langle v', \Y^{f}(e^{(z_{1}-z_{2})L(-1)}\Y^{f}(v, z_{2}-z_{1})
\omega, z_{2})\one\rangle}\nn
&&=\langle e^{(z_{1}-z_{2})L'(-1)}v', \Y^{f}(\Y^{f}(v, z_{2}-z_{1})
\omega, z_{2})\one\rangle
\end{eqnarray}
is absolutely convergent in the region $|z_{2}|>|z_{1}-z_{2}|>0$ and 
is equal to (\ref{iter2-omega}) in the region 
$|z_{1}|, |z_{2}|>|z_{1}-z_{2}|>0$. 
So (\ref{iter-omega}) and the left-hand side of 
(\ref{iter3-omega}) are analytic 
extensions of each other. 

We know that both (\ref{iter-omega}) and the left-hand side 
of (\ref{iter3-omega})
are convergent absolutely in the region $|z_{2}|>|z_{1}-z_{2}|>0$ and,
moreover, we know that 
(\ref{prod1-omega}), (\ref{prod2-omega}), (\ref{iter-omega}) and
(\ref{iter2-omega})
give single-valued analytic functions in $z_{1}$ and $z_{2}$. 
Thus in the region $|z_{2}|>|z_{1}-z_{2}|>0$, 
(\ref{iter-omega}) and the left-hand side of 
(\ref{iter3-omega}) are equal,
that is,
\begin{equation}\label{iter-eqn-1}
\langle v', \Y^{f}(\Y^{f}(\omega, z_{1}-z_{2})v, z_{2})\one\rangle
=\langle v', \Y^{f}(e^{(z_{1}-z_{2})L(-1)}\Y^{f}(v, z_{2}-z_{1})
\omega, z_{2})\one\rangle.
\end{equation}
By taking coefficients of $z_{1}-z_{2}$ and $z_{2}$ 
in both sides of (\ref{iter-eqn-1}) and then taking the generating functions
of these coefficients, we obtain 
\begin{equation}\label{iter-eqn-2}
\langle v', \Y^{f}(\Y^{f}(\omega, x)v, y)\one\rangle
=\langle v', \Y^{f}(e^{xL(-1)}\Y^{f}(v, -x)
\omega, y)\one\rangle,
\end{equation}
where $x$ and $y$ are commuting formal variables.
Since $v'\in V'$ is arbitrary, (\ref{iter-eqn-2})
gives 
\begin{equation}\label{iter-eqn-3}
\Y^{f}(\Y^{f}(\omega, x)v, y)\one
=\Y^{f}(e^{xL(-1)}\Y^{f}(v, -x)
\omega, y)\one.
\end{equation}
Taking the formal limit $y\to 0$
(that is, taking the constant term of 
of series in $y$)
of both sides of (\ref{iter-eqn-3}), we obtain 
$$\Y^{f}(\omega, x)v=e^{xL(-1)}\Y^{f}(v, -x)\omega.$$
So we conclude that $\omega\in C_{0}(V)$.
\epfv

One immediate consequence of this result is the following:

\begin{cor}
Let $V$ be a grading-restricted conformal open-string
vertex algebra. Then the vertex operator algebra
$\langle \omega\rangle$ is a subalgebra of $C_{0}(V)$.
\end{cor}

Recall the following main theorem in \cite{H7}:

\begin{thm}\label{ioa}
Let $V$ be a vertex operator algebra satisfying the 
following conditions:

\begin{enumerate}

\item  Every  generalized $V$-module is a direct sum of 
irreducible $V$-modules.

\item There are only finitely many inequivalent irreducible $V$-modules
and these irreducible $V$-modules are all $\R$-graded.  

\item Every
irreducible $V$-module satisfies the $C_{1}$-cofiniteness condition.

\end{enumerate}

\noindent Then the direct sum of all (inequivalent) irreducible
$V$-modules has a natural structure of an intertwining operator algebra.
In particular,
the following associativity for intertwining operators holds:
For any $V$-modules 
$W_{0}$, $W_{1}$, $W_{2}$, $W_{3}$  and $W_{4}$, 
any 
intertwining operators ${\cal Y}_{1}$ and ${\cal Y}_{2}$ of 
types $\binom{W_{0}}{W_{1}W_{4}}$ and $\binom{W_{4}}{W_{2}W_{3}}$,
respectively,
\begin{equation}\label{int-prod}
\langle w'_{(0)}, {\cal Y}_{1}(w_{(1)},
z_{1}){\cal Y}_{2}(w_{(2)}, z_{2})w_{(3)}\rangle
\end{equation}
is absolutely convergent when $|z_{1}|>|z_{2}|>0$ for
$w'_{(0)}\in W'_{0}$, $w_{(1)}\in W_{1}$,
$w_{(2)}\in W_{2}$ and $w_{(3)}\in W_{3}$, and there exist
$V$-module $W_{5}$ and 
intertwining operators ${\cal Y}_{3}$
and ${\cal Y}_{4}$ of types $\binom{W_{5}}{W_{1}W_{2}}$ and 
$\binom{W_{0}}{W_{5}W_{3}}$, respectively, such that 
$$\langle w'_{(0)}, {\cal Y}_{4}({\cal Y}_{3}(w_{(1)},
z_{1}-z_{2})w_{(2)}, z_{2})w_{(3)}\rangle$$
is absolutely convergent when $|z_{2}|>|z_{1}-z_{2}|>0$
for $w'_{(0)}\in W'_{0}$, $w_{(1)}\in W_{1}$,
$w_{(2)}\in W_{2}$ and $w_{(3)}\in W_{3}$ and is equal to 
(\ref{int-prod}) when $|z_{1}|>|z_{2}|>|z_{1}-z_{2}|>0$.
\end{thm}

Theorems \ref{op-va-int} and \ref{ioa} suggest a 
method to construct conformal open-string vertex algebra: We 
start with a vertex operator algebra $(V, Y, \mathbf{1}, \omega)$ 
satisfying the conditions 
in Theorem \ref{ioa} and look for a module $W$ and 
an intertwining operator $\mathcal{Y}^{f}$ of type $\binom{W}{WW}$ 
such that if we define 
\begin{eqnarray*}
Y^{O}: (W\otimes W)\times \R_{+}&\to& \overline{W}\\
(w_{1}\otimes w_{2}, r)&\mapsto& Y^{O}(w_{1}, r)w_{2}
\end{eqnarray*}
by
\begin{equation}\label{y-r}
Y^{O}(w_{1}, r)w_{2}=
\Y^{f}(w_{1}, r)w_{2}\end{equation}
for $r\in \R_{+}$,
then $(W, Y^{O}, \mathbf{1}, \omega)$ is a conformal 
open-string vertex algebra.

We give more details here. Let $(V, Y, \mathbf{1}, \omega)$ 
be a vertex operator algebra
satisfying the conditions in Theorem \ref{ioa}. For simplicity, we 
assume that $V$ is simple. Let
$\mathcal{A}$ be the
set of equivalence classes of 
irreducible $V$-modules and, for $a\in \mathcal{A}$, let $W^{a}$ be
a representative in $a$. Then by Theorem \ref{ioa},
$\coprod_{a\in A} W^a$ has a natural structure of
an intertwining operator algebra.
Let $W=\coprod_{a\in A} E^a\otimes W^a$ where $E^{a}$ for $a\in \mathcal{A}$
are vector spaces to be determined. We give $W$ the 
obvious $V=\mathbb{C}\otimes V$-module structure. 
We also let
\begin{eqnarray*}
\Y^{f} &\in& \text{Hom}(W\otimes W, W\{ x\})\nn
&&=\coprod_{a_1, a_2, a_3\in \mathcal{A}}
\text{Hom}(E^{a_{1}}\otimes E^{a_{2}}, 
E^{a_{3}})\otimes \text{Hom}(W^{a_{1}}\otimes W^{a_{2}}, 
W^{a_{3}}\{ x\})
\end{eqnarray*}
be given by
$$\Y^{f} = \sum_{a_1,a_2,a_3\in \mathcal{A}} 
\sum_{i=1}^{\mathcal{N}_{a_{1}a_{2}}^{a_{3}}}
C_{a_1a_2}^{a_3;i} 
\otimes \Y_{a_1a_2}^{a_3;i}
$$ 
where  for $a_{1}, a_{2}, a_{3}\in \mathcal{A}$, 
$\mathcal{N}_{a_{1}a_{2}}^{a_{3}}$
is the fusion rule
of type $\binom{W^{a_{3}}}{W^{a_{1}}W^{a_{2}}}$,
$C_{a_1a_2}^{a_3;i}\in \text{Hom}(E^{a_{1}}\otimes E^{a_{2}}, 
E^{a_{3}})$ for $i=1, \dots, \mathcal{N}_{a_{1}a_{2}}^{a_{3}}$ 
are to be determined, and 
$\Y_{a_1a_2}^{a_3;i}$ for $i=1, \dots, \mathcal{N}_{a_{1}a_{2}}^{a_{3}}$
is a basis of the space $\mathcal{V}_{a_1a_2}^{a_3}$ of 
intertwining operators of type $\binom{W^{a_{3}}}{W^{a_{1}}W^{a_{2}}}$.

Let $e$ be the equivalence class of irreducible $V$-modules containing 
$V$. Note that $\mathcal{N}_{ea}^{a}$ for $a\in \mathcal{A}$ are always 
one-dimensional.
We choose the basis $\Y_{ea}^{a; 1}$ for  $a\in \mathcal{A}$
to be the vertex operator for the $V$-module $W^{a}$. In particular,
$\Y_{ea}^{a; 1}(\one, x)w^{a}=w^{a}$ for $a\in \mathcal{A}$ and 
$w^{a}\in W^{a}$. We also choose the basis 
$\Y_{ae}^{a; 1}$ for  $a\in \mathcal{A}$
to be the ones given by 
$$\Y_{ae}^{a; 1}(w^{a}, x)u=e^{xL(-1)}\Y_{ea}^{a; 1}(u, -x)w$$
for $u\in V$ and $w^{a}\in W^{a}$. Thus we have
$\lim_{x\to 0}\Y_{ae}^{a; 1}(w^{a}, x)\one=w^{a}$ 
for $a\in \mathcal{A}$ and 
$w^{a}\in W^{a}$.

We would like to choose $E^{a}$ for $a\in \mathcal{A}$
and  $C_{a_1a_2}^{a_3; i}$ for $a_{1}, a_{2}, a_{3}\in \mathcal{A}$
and $i=1, \dots, \mathcal{N}_{a_{1}a_{2}}^{a_{3}}$  such that 
the map $Y^{O}$ given by (\ref{y-r}) in terms of 
$\Y^{f}$ satisfies the
associativity 
\begin{equation}\label{assoc-eqn}
Y^{O}(w_1, r_1)Y^{O}(w_2, r_2)w_3 = Y^{O}(Y^{O}(w_1, r_1-r_2)w_2, r_2)w_3
\end{equation}
for $r_{1}, r_{2}\in \R_{+}$ satisfying $r_{1}>r_{2}>r_{1}-r_{2}>0$
and $w_{1}\in W_{1}$, $w_{2}\in W_{2}$, $w_{3}\in W_{3}$. Note that 
both sides of (\ref{assoc-eqn}) are well-defined 
since $\coprod_{a\in A} W^a$ is an intertwining operator 
algebra. 
The left-hand side of (\ref{assoc-eqn})
gives 
\bea
\lefteqn{\sum_{\substack{a_1,a_2,a_3}{a_4,a;i,j}}
(C_{a_1a}^{a_4;i}\circ (\id_{E^{a_{1}}}\otimes C_{a_2a_3}^{a;j}))
\otimes \Y_{a_1a}^{a_4;i}(w_1,r_1)\Y_{a_2a_3}^{a;j}(w_2,r_2)w_3}   \nn
&&= \sum_{\substack{a_1,a_2,a_3}{a_4,a;i,j}}
(C_{a_1a}^{a_4;i}\circ (\id_{E^{a_{1}}}\otimes C_{a_2a_3}^{a;j}))\nn
&&\quad\quad\quad\quad\quad\otimes \sum_{a_{5};k,l} 
\mathcal{F}_{a;a_{5}}^{ij;kl}(a_1,a_2,a_3;a_4)
\Y_{a_{5}a_3}^{a_4;l}(\Y_{a_1a_2}^{a_{5};k}(w_1,r_1-r_2)w_2,r_2)w_3  \nonumber
\eea
where for any $a\in \mathcal{A}$, $\id_{E^{a}}$ is the identity on 
$E^{a}$ and $\mathcal{F}_{a;a_{5}}^{ij;kl}(a_1,a_2,a_3;a_4)$,
for $a, a_{1}, \dots, a_{5}\in \mathcal{A}$, 
$i=1, \dots, \mathcal{N}_{a_{1}a}^{a_{4}}$,
$j=1, \dots, \mathcal{N}_{a_{2}a_{3}}^{a}$,
$k=1, \dots, \mathcal{N}_{a_{1}a_{2}}^{a_{5}}$ and
$l=1, \dots, \mathcal{N}_{a_{5}a_{3}}^{a_{4}}$,
are the matrix elements of the corresponding fusing isomorphisms.
(In the formulas above and below, for simplicity, we omit
the ranges over which the sums are taken, since these are clear and
some of them have been given above.)
The right-hand side of (\ref{assoc-eqn}) gives
$$\sum_{\substack{a_1,a_2,a_3}{a_4,a_{5};k,l}}
(C_{a_{5}a_3}^{a_4;l}\circ (C_{a_1a_2}^{a_{5};k}\otimes \id_{E^{a_{3}}}))
\otimes \Y_{a_{5}a_3}^{a_4;l}(\Y_{a_1a_2}^{a_{5};k}(w_1,r_1-r_2)w_2,r_2)w_3.$$
It is clear that in this case
$\Y_{a_{5}a_3}^{a_4;l}(\Y_{a_1a_2}^{a;k}
(\cdot,r_1-r_2)\cdot,r_2)\;\cdot$ for $a_1, a_2, a_3, a_4, a_{5}
\in \mathcal{A}$
are linearly independent. Thus (\ref{assoc-eqn}) gives
\begin{equation}  \label{fusing-1}
\sum_{a;i,j} \mathcal{F}_{a;a_{5}}^{ij;kl}(a_1,a_2,a_3;a_4)
(C_{a_1a}^{a_4;i}\circ (\id_{E^{a_{1}}}\otimes C_{a_2a_3}^{a;j}))
= C_{a_{5}a_3}^{a_4;l}\circ (C_{a_1a_2}^{a_{5};k}\otimes \id_{E^{a_{3}}})
\end{equation}
for $a_{1}, a_{2}, a_{3}, a_{4}, a_{5}\in \mathcal{A}$, 
$k=1, \dots, \mathcal{N}_{a_{1}a_{2}}^{a_{5}}$ and
$l=1, \dots, \mathcal{N}_{a_{5}a_{3}}^{a_{4}}$.

We need a vacuum for $W$. 
Let $1^{e}\in E^{e}$. If we want the vacuum to be of the form 
$\one_{W}=1^{e}\otimes \one$, then we must have the following 
identity property
and creation property:
\begin{eqnarray}
Y^{O}(\one_{W}, r)(\alpha^{a}\otimes w^{a})&=&\alpha^{a}\otimes w^{a},
\label{identity1}\\
\lim_{r\to 0}Y^{O}((\alpha^{a}\otimes w^{a}), r)\one_{W}
&=&\alpha^{a}\otimes w^{a}\label{identity2}
\end{eqnarray}
for $a\in \mathcal{A}$, $\alpha^{a}\in E^{a}$ and $w^{a}\in W^{a}$.
The equations (\ref{identity1}) and (\ref{identity2})
together with the properties of 
intertwining operators for $V$ gives
\begin{eqnarray}
C_{ea}^{a; 1}(1^{e}\otimes \alpha^{a})&=&\alpha^{a},\label{identity3}\\
C_{ae}^{a; 1}(\alpha^{a} \otimes 1^{e})&=&\alpha^{a}\label{identity4}
\end{eqnarray}
for $a\in \mathcal{A}$ and $\alpha^{a}\in E^{a}$.

Let $\mathbf{1}_{W}=1^{e}\otimes \mathbf{1}$ and
$\omega_W=1^{e}\otimes \omega$. Then
we have just proved the following:

\begin{prop}
Let $V$ be a simple vertex operator algebra satisfying the conditions 
in Theorem \ref{ioa} and let $\mathcal{A}$, $e$ and 
$W^{a}$ for $a\in \mathcal{A}$ be as above. 
If we choose the vector spaces $E^{a}$ for $a\in \mathcal{A}$, 
$C_{a_1a_2}^{a_3; i}\in \text{Hom}(E^{a_{1}}\otimes E^{a_{2}}, 
E^{a_{3}})$ for $a_{1}, a_{2}, a_{3}\in \mathcal{A}$,
$i=1, \dots, \mathcal{N}_{a_{1}a_{2}}^{a_{3}}$, and 
$1^{e}\in E^{e}$  such that 
(\ref{fusing-1}), (\ref{identity3}) and (\ref{identity4}) hold, 
then the quadruple $(W,Y^{O},\mathbf{1}_{W},\omega_{W})$ 
is a grading-restricted
conformal open-string vertex algebra. 
\end{prop}

\renewcommand{\theequation}{\thesection.\arabic{equation}}
\renewcommand{\thethm}{\thesection.\arabic{thm}}
\setcounter{equation}{0}
\setcounter{thm}{0}

\section{Examples}

In this section, we give some examples of open-string vertex algebras. 
Examples can also be constructed using the main results in Sections 4
and 5.

First of all, we have the following examples for which the axioms are
trivial to verify:

\begin{enumerate}

\item Associative algebras.

\item Vertex (super)algebras.

\item Tensor products of algebras above, for example,
$A\otimes V$ where $A$ is an associative algebra and $V$ a 
vertex (super)algebra.

\end{enumerate}

The examples above are trivial to construct because they satisfy some
much stronger axioms than those in the definition of open-string vertex
algebra. Nontrivial examples of open-string vertex algebras can be
constructed {}from the direct sum of a vertex algebra and an $\R$-graded
module for the vertex algebra in the same
ways as in the construction of the example of vertex operator algebras in
Example 3.4 in \cite{H3} and as in the conceptual construction of the vertex
operator algebra structure on the moonshine module in \cite{H4.5},
except that here the module does not have to be $\Z$-graded. Note that in
the construction of the vertex operator algebra structure on the
moonshine module in \cite{H4.5}, the hard part is to prove the duality
properties, which follow {}from the duality properties of a larger
intertwining operator algebra. If we start with a vertex operator
algebra satisfying the conditions in Theorem \ref{ioa}, then the
construction becomes very easy because the duality properties have been 
established by Theorem \ref{ioa}. 

We now give an example constructed using a different method. It
is an example constructed {}from modules for the minimal Virasoro vertex
operator algebra of central charge $c=\frac{1}{2}$. This example is
nontrivial because it is not an associative algebra, a vertex
(super)algebra or a tensor product of these algebras.  Here we describe
the data. For the details, we refer the reader to the second author's
thesis \cite{K}. For the minimal Virasoro vertex operator algebras,
their representations, intertwining operators and chiral correlation
functions, see, for example, \cite{DF}, \cite{BPZ}, \cite{W}, \cite{H4},
\cite{FRW} and \cite{DMS}.

Let $L(\frac{1}{2}, 0)$ be the minimal Virasoro vertex operator algebra
of central charge $\frac{1}{2}$. It
has three inequivalent irreducible modules
$W_{0}=L(\frac{1}{2}, 0)$, $W_{1}=L(\frac{1}{2}, \frac{1}{2})$ 
and $W_{2}=L(\frac{1}{2}, \frac{1}{16})$. It is well known that 
the fusion rules $\mathcal{N}_{ij}^{k}=\mathcal{N}_{W_{i}W_{j}}^{W_{k}}$ 
for $i, j, k=0, 1, 2$ are equal to $1$ for
\begin{eqnarray*}
(i, j, k)&=&(0, 0, 0), (0, 1, 1), (1, 0, 1), (1, 1, 0),
(0, 2, 2),\nn
&& (2, 0, 2), (2, 2, 0), (1, 2, 2), (2, 1, 2), (2, 2, 1)
\end{eqnarray*}
and are equal to $0$ otherwise. It was proved in \cite{H4}
that the direct sum of $W_{0}$, $W_{1}$ and $W_{3}$ has a structure of 
intertwining operator algebra. When $\mathcal{N}_{ij}^{k}=1$,
we choose a basis $\Y_{ij}^{k}$ of $\mathcal{V}_{ij}^{k}$.
Given $i, j, k, l\in \{0, 1, 2\}$, $m\in \{0, 1, 2\}$ 
is said to be {\it coupled with $n\in \{0, 1, 2\}$ through}
$(i, j, k; l)$ if 
$\mathcal{V}_{im}^{l}$, 
$\mathcal{V}_{jk}^{m}$, $\mathcal{V}_{ij}^{n}$ and 
$\mathcal{V}_{nk}^{l}$ are all nonzero. We use the notation
$m \Join_{i, j, k}^{l} n$ to denote the fact that $m$ is coupled with $n$
through $(i, j, k; l)$.

For $i, j, k, l\in \{0, 1, 2\}$,
the matrix elements $\mathcal{F}_{m; n}(i, j, k; l)$ for
$m, n=0, 1, 2$ 
of the fusing isomorphisms
$$\mathcal{F}(i, j, k; l): \coprod_{m=0}^{2}
\mathcal{V}_{im}^{l}\otimes 
\mathcal{V}_{jk}^{m}\to \coprod_{n=0}^{2}\mathcal{V}_{ij}^{n}\otimes 
\mathcal{V}_{nk}^{l}$$
are determined by the following associativity relations (see \cite{H6})
\begin{eqnarray*}
\lefteqn{\langle w'_{l}, \Y_{im}^{l}(w_{i},
z_{1})\Y_{jk}^{m}(w_{j}, z_{2})w_{k}\rangle}\nn
&&=\sum_{m\Join_{i, j, k}^{l} n}\mathcal{F}_{m; n}(i, j, k; l)
\langle w'_{l}, \Y_{nk}^{l}(\Y_{ij}^{n}(w_{i}, z_{1}-z_{2})
w_{j}, z_{2})w_{k}\rangle
\end{eqnarray*}
for $i, j, m=0, 1, 2$,
$z_{1}, z_{2}\in \R$ satisfying $z_{1}>z_{2}>z_{1}-z_{2}>0$
and $w_{i}\in W_{i}$, $w_{j}\in W_{j}$, $w_{k}\in W_{k}$,
where the sum is over all $k, l, n=0, 1, 2$ 
such that $m\Join_{i, j, k}^{l} n$. For simplicity, we use
$\tilde{\mathcal{F}}(i, j, k; l)$ for $i, j, k, l=0, 1, 2$
to denote  matrices whose entries 
$\tilde{\mathcal{F}}_{mn}(i, j, k; l)$ for 
$m, n=0, 1, 2$ is the symbol $DC$ (meaning decoupled)
if $m$ is not 
coupled with $n$ through $(i, j, k; l)$
and is $\mathcal{F}_{m; n}(i, j, k; l)$
if $m$ is coupled with $n$ through $(i, j, k; l)$. We call these matrices the 
{\it fusing-coupling} matrices.
For $m, n=0, 1, 2$, 
we use $\pm E_{mn}$ to denote the $3\times 3$  matrices 
with the entry in the $m$-th row and the $n$-th column being $\pm 1$ 
and the other entries being $DC$. 

\begin{prop}
For $i, j, k=0, 1, 2$ such that $\mathcal{N}_{ij}^{k}=1$,
there exist basis $\Y_{ij}^{k}$ of  
$\mathcal{V}_{ij}^{k}$
such that 
\begin{eqnarray*}
&\tilde{\mathcal{F}}(0, 0, 0, 0)=\tilde{\mathcal{F}}(1, 1, 1, 1)=E_{00},&\\
&\tilde{\mathcal{F}}(1, 1, 0, 0)=\tilde{\mathcal{F}}(0, 0, 1, 1)=E_{01},&\\
&\tilde{\mathcal{F}}(1, 0, 0, 1)=\tilde{\mathcal{F}}(0, 1, 1, 0)=E_{10},&\\
&\tilde{\mathcal{F}}(1, 0, 1, 0)=\tilde{\mathcal{F}}(0, 1, 0, 1)=E_{11},&\\
&\tilde{\mathcal{F}}(2, 2, 0, 0)=\tilde{\mathcal{F}}(0, 0, 2, 2)
=\tilde{\mathcal{F}}(1, 1, 2, 2)=\tilde{\mathcal{F}}(2, 2, 1, 1)=E_{02},&\\
&\tilde{\mathcal{F}}(0, 2, 2, 0)=\tilde{\mathcal{F}}(2, 0, 0, 2)
=\tilde{\mathcal{F}}(1, 2, 2, 1)=\tilde{\mathcal{F}}(2, 1, 1, 2)=E_{20},&\\
&\tilde{\mathcal{F}}(0, 1, 2, 2)=\tilde{\mathcal{F}}(1, 0, 2, 2)
=\tilde{\mathcal{F}}(2, 2, 0, 1)=\tilde{\mathcal{F}}(2, 2, 1, 0)=E_{12},&\\
&\tilde{\mathcal{F}}(0, 2, 2, 1)=\tilde{\mathcal{F}}(1, 2, 2, 0)
=\tilde{\mathcal{F}}(2, 0, 1, 2)=\tilde{\mathcal{F}}(2, 1, 0, 2)=E_{21},&\\
&\quad\tilde{\mathcal{F}}(1, 2, 1, 2)=\tilde{\mathcal{F}}(2, 1, 2, 1)
=-E_{22},&\\
&\tilde{\mathcal{F}}(0, 2, 0, 2)=\tilde{\mathcal{F}}(2, 0, 2, 0)
=\tilde{\mathcal{F}}(0, 2, 1, 2)=\tilde{\mathcal{F}}(1, 2, 0, 2)
\quad\quad\quad&\\
&\quad\quad\quad\quad\quad\quad\quad\quad\quad\quad\quad\;\;
=\tilde{\mathcal{F}}(2, 0, 2, 1)=\tilde{\mathcal{F}}(2, 1, 2, 0)=E_{22},\\
&\;\;\tilde{\mathcal{F}}(2, 2, 2, 2)=\left(\begin{array}{ccc}
\frac{1}{\sqrt{2}}&\frac{1}{\sqrt{2}}&DC\\
\frac{1}{\sqrt{2}}&-\frac{1}{\sqrt{2}}&DC\\
DC&DC&DC
\end{array}\right);&\\
\end{eqnarray*}
all other fusing-coupling matrices have entries 
which are either $0$ or $DC$.
\end{prop}

The proposition above gives the  complete information about the 
fusing isomorphisms for the minimal model of central charge $\frac{1}{2}$.
Now consider the irreducible modules 
$W_i\otimes W_i$ for $i=0, 1, 2$ for the 
tensor product vertex operator algebra $L(\frac{1}{2}, 0)\otimes 
L(\frac{1}{2}, 0)$.
Let $W=\coprod_{i=0}^2 W_i\otimes W_i$ and let
$$\Y^{f}: (W\otimes W) \to W\{ x\}$$
be given by
$$\Y^{f}= \sum_{i, j, k=0}^2 \Y_{ij}^k\otimes \Y_{ij}^k$$
where we have taken $\Y_{ij}^{k}=0$
for $i, j, k\in \{0, 1, 2\}$ such that $\mathcal{V}_{ij}^{k}=0$ and 
where $\Y_{ij}^k\otimes \Y_{ij}^k$ for $i, j, k\in \{0, 1, 2\}$
act on $W\otimes W$ in the obvious way. 
Let 
\begin{eqnarray*}
Y^{O}: (W\otimes W)\times \R_{+} &\to& \overline{W}\nn
(w_{1}\otimes w_{2}, r)&\mapsto& Y^{O}(w_{1}, r)w_{2}
\end{eqnarray*}
be given by 
$$Y^{O}(w_{1}, r)w_{2}=\Y^{f}(w_{1}, r)w_{2}$$
for $r\in \R_{+}$ and $w_{1}, w_{2}\in W$.
Let $\one$ and $\omega$ be the vacuum and conformal element
of $L(\frac{1}{2}, 0)$. Then we have:

\begin{prop}
The quadruple $(W,Y^{O},\mathbf{1}\otimes \mathbf{1}, 
\omega\otimes \mathbf{1}+\mathbf{1}\otimes 
\omega)$ is a grading-restricted conformal open-string 
vertex algebra with
$C_0(W)=W_0\otimes W_0$. 
\end{prop}
 
The proof is a straightforward verification. See \cite{K} for 
details.

\begin{rema}
{\rm In the construction above, $\Y^{f}$ and 
$Y^{O}$ involve fractional powers. So
$W$ is not
a vertex operator algebra.}
\end{rema}

\renewcommand{\theequation}{\thesection.\arabic{equation}}
\renewcommand{\thethm}{\thesection.\arabic{thm}}
\setcounter{equation}{0}
\setcounter{thm}{0}

\section{Braided tensor categories and open-string vertex algebras}

In this section, we show that an associative algebra in 
the braided tensor category of modules for a suitable 
vertex operator algebra $V$
is equivalent to an open-string vertex algebra with 
$V$ in its meromorphic center. 
The main result of this section (Theorem \ref{main}) is a
straightforward generalization of the main result in 
\cite{HKL}. 
In this section, we assume that the reader is familiar 
with the tensor product theory developed by Lepowsky and the first author.
See \cite{HL3}--\cite{HL6} and \cite{H3} for details.

First of all, we have the following result established in \cite{H8}:

\begin{thm}\label{tensor-cat}
Let $V$ be a vertex operator algebra satisfying the 
conditions in Theorem \ref{ioa}. Then the 
category of $V$-modules has a natural structure of 
vertex tensor category with $V$ as its unit object. 
In particular, this category has a natural 
structure of braided tensor category.
\end{thm}

Given a braided tensor category $\mathcal{C}$, we use $1_{\mathcal{C}}$ 
to denote its 
unit object. We need the following concept:

\begin{defn}\label{c-alg}
{\rm Let $\mathcal{C}$ be a braided tensor category. 
An {\it associative algebra in $\mathcal{C}$} (or 
{\it associative $\mathcal{C}$-algebra}) 
is an object $A\in \mathcal{C}$ along with a morphism
$\mu: A\otimes
A\to A$ and an injective morphism 
$\iota_A: 1_{\mathcal{C}} \to A$ such that the following
conditions hold:

\begin{enumerate}

\item Associativity:
$\mu \circ(\mu \otimes \id_{A})=\mu\circ 
(\id_{A}\otimes \mu)\circ \mathcal{A}$
where $\mathcal{A}$ is the associativity isomorphism {}from
$A\otimes (A\otimes A)$ to $(A\otimes A)\otimes A$.

\item Unit properties:
$\mu \circ (\iota_A \otimes \id_A)\circ l^{-1}_{A}
=\mu \circ (\iota_{A}\otimes  \id_{A})\circ r^{-1}_{A}=\id_A$
where $l_{A}: 1_{\mathcal{C}}\otimes A\to A$ and $r_{A}: A\otimes 
1_{\mathcal{C}}\to A$
are the left and right unit isomorphism, respectively.

\end{enumerate}

We say that the unit of 
an associative algebra $A$ in $\mathcal{C}$ is {\it unique}
if 
$$\dim \hom_{\mathcal{C}}(1_{\mathcal{C}}, A)=1.$$ 
}
\end{defn}

We use $(A, \mu, \iota_{A})$ or simply $A$ to denote 
the associative algebra in $\mathcal{C}$ just defined.

Let $V$ be a vertex operator algebra satisfying 
the conditions in Theorem \ref{ioa}. Then we know that the 
direct sum of all irreducible $V$-modules is an intertwining operator
algebra. 
We say that this intertwining operator algebra  satisfies the
{\it positive weight condition} if
for any irreducible $V$-module $W$, the weights of nonzero elements
of $W$ are nonnegative, $W_{(0)}\ne 0$ if and only if 
$W$ is isomorphic to $V$, 
and $V_{(0)}=\mathbb{C}\mathbf{1}$.
We say that an open-string vertex algebra $V$
satisfy the {\it positive weight condition} if
the weights of elements
of $V$ are nonnegative and
$V_{(0)}=\mathbb{C}\mathbf{1}$.

\begin{thm}\label{main}
Let $(V, Y, \one, \omega)$ be a vertex operator algebra satisfying 
the conditions in Theorem \ref{ioa}
and let $\mathcal{C}$ be the braided tensor category of 
$V$-modules. Then the  categories of the following objects are
isomorphic: 

\begin{enumerate}
\item A grading-restricted conformal open-string vertex algebra $V_{e}$
and an injective homomorphism of vertex operator algebras
{}from $V$ to the meromorphic center $C_{0}(V_{e})$ 
of $V_{e}$.
 
\item An associative algebra $V_{e}$  in $\mathcal{C}$. 
\end{enumerate}

If the intertwining operator algebra on the 
direct sum of all irreducible $V$-modules satisfies the
positive weight condition, then an algebra $V_e$ in Category 1 above 
satisfy the positive weight condition if and only 
if the unit of the corresponding associative algebra $V_{e}$ 
in $\mathcal{C}$ is unique.
\end{thm}
\pf
Let $V_e$ be a grading-restricted conformal 
open-string vertex  algebra, $\one_{e}$ the vacuum of $V_{e}$
and $\iota_{V_{e}}$ 
an injective homomorphism of 
vertex operator algebras {}from $V$ to $C_{0}(V_{e})$. Then 
we have $\iota_{V_{e}}(\one)=\one_{e}$.
Then by Theorem \ref{op-va-int}, 
$V_e$ is an $\iota_{V_{e}}(V)$-module
and thus a $V$-module. So $V_e$ is an
object in $\mathcal{C}$. Since $V_e$ is an open-string vertex  algebra, 
we have a vertex operator map $Y^{O}_e$ for $V_e$. By 
Theorem \ref{op-va-int} again, the corresponding 
formal vertex operator map 
$\Y^{f}_{e}$ is in fact an intertwining operator for $V$ of type 
$\binom{V_{e}}{V_{e}V_{e}}$.   Let 
$\mu: V_{e}\boxtimes V_{e}\to V_{e}$ be the module map corresponding to 
the intertwining operator $\Y^{f}_{e}$. We claim that
$(V_{e}, \mu,  \iota_{V_{e}})$ is  an associative algebra 
in $\mathcal{C}$. The proof is similar to
the proof of the result in \cite{HKL} that 
suitable commutative associative algebras
in $\mathcal{C}$ are equivalent to
vertex operator algebras extending $V$. For reader's convenience,
we give a proof  here.

For $r\in \R_{+}$, 
let $\mu_{r}$ be the morphism {}from 
$V_{e}\boxtimes_{P(r)}V_{e}$ to $V_{e}$ corresponding to the intertwining 
operator $Y^{f}_{e}$ and let
$\overline{\mu}_{r}: \overline{V_{e}\boxtimes_{P(r)} 
V_{e}}\to \overline{V}_{e}$
be the natural extension of $\mu_{P(r)}$.
Then by definition, $\mu=\mu_{1}$ and
$$\overline{\mu}_{r}(u\boxtimes_{P(r)} v)
=\Y_{e}^{f}(u, r)v=Y^{O}_{e}(u, r)v.$$
for $u, v\in V_{e}$.  For simplicity, we shall use $\id$ to denote
$\id_{V_{e}}$ in this proof. 
Thus for $u, v, w\in V_{e}$ and $r_{1}, r_{2}\in \R_{+}$ satisfying
$r_{1}>r_{2}>r_{1}-r_{2}>0$, 
\begin{eqnarray}
\lefteqn{\overline{(\mu_{r_{1}} \circ(\id \boxtimes_{P(r_{1})} \mu_{r_{2}}))}
(u\boxtimes_{P(r_{1})}
(v\boxtimes_{P(r_{2})} w))}\nn
&&\quad\quad\quad\quad =Y^{O}_{e}(u, r_{1})Y^{O}_{e}(v, r_{2})w
\label{prod-tensor}\\
\lefteqn{\overline{(\mu_{r_{2}} \circ(\mu_{r_{1}-r_{2}} \boxtimes_{P(r_{2})} \id))}
((u\boxtimes_{P(r_{1}-r_{2})}
v)\boxtimes_{P(r_{2})} w)}\nn
&&\quad\quad\quad\quad=Y^{O}_{e}Y^{O}_{e}(u, r_{1}-r_{2})v, r_{2})w,
\label{iter-tensor}
\end{eqnarray}
where (and below) we use the notation that 
a linear map preserving gradings with a horizontal 
line over it always mean 
the natural extension of the map to a map between the algebraic completions 
of the original graded spaces.
The associativity for 
$Y^{O}_{e}$ gives 
\begin{equation}
Y^{O}_{e}(u, r_{1})Y^{O}_{e}(v, r_{2})w
=Y^{O}_{e}(Y^{O}_{e}(u, r_{1}-r_{2})v, r_{2})w.\label{assoc}
\end{equation}
The associativity isomorphism 
$$\mathcal{A}^{P(r_{1}-r_{2}), P(r_{2})}_{P(r_{1}), P(r_{2})}: 
V_{e}\boxtimes_{P(r_{1})} (V_{e}\boxtimes_{P(r_{2})} V_{e})\to 
(V_{e}\boxtimes_{P(r_{1}-r_{2})} V_{e})\boxtimes_{P(r_{2})} V_{e}$$
is characterized by 
\begin{equation}
\overline{\mathcal{A}}^{P(r_{1}-r_{2}), P(r_{2})}_{P(r_{1}), P(r_{2})}
(u\boxtimes_{P(r_{1})} (v\boxtimes_{P(r_{2})} w))
=(u\boxtimes_{P(r_{1}-r_{2})} v)\boxtimes_{P(r_{2})} w\label{assoc-iso}
\end{equation}
for $u,v, w\in V_{e}$, where 
$\overline{\mathcal{A}}^{P(r_{1}-r_{2}), P(r_{2})}_{P(r_{1}), P(r_{2})}$
is the natural extension 
of $\mathcal{A}^{P(r_{1}-r_{2}), P(r_{2})}_{P(r_{1}), P(r_{2})}$.

Combining (\ref{prod-tensor})--(\ref{assoc-iso}), we obtain
\begin{equation}
(\mu_{r_{1}} \circ(\id \boxtimes_{P(r_{1})} \mu_{r_{2}}))=
(\mu_{r_{2}} \circ(\mu_{r_{1}-r_{2}} \boxtimes_{P(r_{2})} \id))\circ
\mathcal{A}^{P(r_{1}-r_{2}), P(r_{2})}_{P(r_{1}), P(r_{2})}.\label{mu-assoc1}
\end{equation}
{}From (\ref{mu-assoc1}), we obtain
\begin{eqnarray}
\lefteqn{(\mu_{r_{1}} \circ(\id \boxtimes_{P(r_{1})} \mu_{r_{2}}))
\circ (\id \boxtimes_{P(r_{1})} 
\mathcal{T}_{\gamma_{2}})\circ \mathcal{T}_{\gamma_{1}}}\nn
&&=
(\mu_{r_{2}} \circ(\mu_{r_{1}-r_{2}} \boxtimes_{P(r_{2})} \id))\circ
\mathcal{A}^{P(r_{1}-r_{2}), P(r_{2})}_{P(r_{1}), P(r_{2})}
\circ (\id \boxtimes_{P(r_{1})} 
\mathcal{T}_{\gamma_{2}})\circ \mathcal{T}_{\gamma_{1}}.\label{mu-assoc2}\nn
&&
\end{eqnarray}
where $r_{1}, r_{2}$ are real numbers satisfying 
$r_{1}>r_{2}>r_{1}-r_{2}> 0$, $\gamma_{1}$ and $\gamma_{2}$ are 
paths in $\R_{+}$
{}from $1$ to $r_{1}$ and $r_{2}$, respectively, 
and $\mathcal{T}_{\gamma_{1}}$ and $\mathcal{T}_{\gamma_{2}}$ the 
parallel transport isomorphisms associated to $\gamma_{1}$ and $\gamma_{2}$, 
respectively. (For reader's convenience, we recall the 
definition of parallel transport isomorphism here. Let 
$\gamma$ be a path {}from $z_{1}\in \C^{\times}$ to $z_{2}\in \C^{\times}$.
The parallel isomorphism 
$\mathcal{T}_{\gamma}: W_{1}\boxtimes_{P(z_{1})}W_{2}\to 
W_{1}\boxtimes_{P(z_{2})}W_{2}$ is given as follows:
Let $\mathcal{Y}$ be the intertwining 
operator corresponding to the intertwining map $\boxtimes_{P(z_{2})}$
and $l(z_{1})$ the value of the logarithm of $z_{1}$ determined uniquely by
$\log z_{2}$ (satisfying $0\le \Im(\log z_{2})<2\pi$) and the 
path $\gamma$. Then  $\mathcal{T}_{\gamma}$ is characterized by
$$\overline{\mathcal{T}}_{\gamma}(w_{1}
\boxtimes_{P(z_{1})}w_{2})=\mathcal{Y}(w_{1}, x)w_{2}
\lbar_{x^{n}=e^{nl(z_{1})}, \;n\in \C}$$
for $w_{1}\in W_{1}$ and $w_{2}\in W_{2}$, where 
$\overline{\mathcal{T}}_{\gamma}$
is the natural extension of $\mathcal{T}_{\gamma}$ 
to the algebraic 
completion $\overline{W_{1}\boxtimes_{P(z_{1})} W_{2}}$ 
of $W_{1}\boxtimes_{P(z_{1})} W_{2}$.  The parallel isomorphism 
depends only on the homotopy class of $\gamma$.)

By definition, we have
\begin{equation}
(\mu_{r_{1}} \circ(\id \boxtimes_{P(r_{1})} \mu_{r_{2}}))\circ 
(\id \boxtimes_{P(r_{1})} 
\mathcal{T}_{\gamma_{2}})\circ \mathcal{T}_{\gamma_{1}}
=\mu \circ(\id \boxtimes \mu).\label{mu-assoc3}
\end{equation}
Similarly, we have
\begin{equation}
(\mu_{r_{2}} \circ(\mu_{r_{1}-r_{2}} \boxtimes_{P(r_{2})} \id))\circ
(\mathcal{T}_{\gamma_{3}}\circ (\mathcal{T}_{\gamma_{4}}
\boxtimes_{P(r_{2})} \id))^{-1}
=(\mu \circ(\mu \boxtimes \id)),\label{mu-assoc4}
\end{equation}
where $\gamma_{3}$ and $\gamma_{4}$ are 
paths in $\R_{+}$
{}from  $r_{2}$ and $r_{1}-r_{2}$ to $1$, respectively, 
and $\mathcal{T}_{\gamma_{3}}$ and $\mathcal{T}_{\gamma_{4}}$ the 
parallel transport isomorphisms associated to $\gamma_{3}$ and $\gamma_{4}$, 
respectively. 
Combining (\ref{mu-assoc2})--(\ref{mu-assoc4}) with the definition
\begin{equation}
\mathcal{A}=\mathcal{T}_{\gamma_{3}}\circ (\mathcal{T}_{\gamma_{4}}
\boxtimes_{P(z_{2})} \id)\circ 
\mathcal{A}^{P(z_{1}-z_{2}), P(z_{2})}_{P(z_{1}), P(z_{2})}\circ
(\id \boxtimes_{P(z_{1})} 
\mathcal{T}_{\gamma_{2}})\circ \mathcal{T}_{\gamma_{1}}.\label{mu-assoc5}
\end{equation}
of the associativity isomorphism for the tensor product structure, 
we obtain the associativity
$$\mu \circ(\id \boxtimes \mu)=(\mu \circ(\mu \boxtimes \id))\circ 
\mathcal{A}.$$

For the unit property, we note that the inverse 
$l_{V_{e}}^{-1}: V_{e}\to V\boxtimes V_{e}$ of the left unit isomorphism 
is defined by
$l_{V_{e}}^{-1}(u)=\mathbf{1}\boxtimes u$ for $u\in V_{e}$ and thus
\begin{eqnarray*}
(\mu\circ (\iota_{V_{e}}\boxtimes \id_{V_{e}})\circ l_{V_{e}}^{-1})(u)
&=&\mu((\iota_{V_{e}}\boxtimes \id_{V_{e}})(\mathbf{1}\boxtimes u))\nn
&=&\mu(\mathbf{1}_{e}\boxtimes u)\nn
&=&Y_{e}(\mathbf{1}_{e}, 1)u\nn
&=&\id_{V_{e}}(u)
\end{eqnarray*}
for $u\in V_{e}$. The other unit property is proved similarly.

Conversely, 
let $(V_{e}, \mu, \iota_{V_{e}})$ be an associative $\mathcal{C}$-algebra. 
In particular,
$V_{e}$ is a $V$-module. The 
module map $\mu: V_{e}\boxtimes V_{e}\to V_{e}$ corresponds to an
intertwining operator $\Y^{f}_{e}$ of type $\binom{V_{e}}{V_{e}V_{e}}$
such that 
\begin{equation}
\overline{\mu}(u\boxtimes v)=\Y^{f}_{e}(u, 1)v\label{converse1}
\end{equation}
for $u, v\in V_{e}$.
Let $\one_{e}=\iota_{V_{e}}(\one)$ and 
$\omega_{e}=\iota_{V_{e}}(\omega)$. We define 
\begin{eqnarray*}
Y^{O}_{e}: (V_{e}\otimes V_{e})\times \R_{+}
&\to& \overline{V}_{e}\\
(u\otimes v, r)&\mapsto&  Y^{O}_{e}(u, r)v
\end{eqnarray*}
by 
$$Y^{O}_{e}(u, r)v=\Y^{f}_{e}(u, r)v$$
for $r\in \R_{+}$, $u, v\in V_{e}$.
Then we claim that
$(V_e, Y^{O}_e,\one_{e}, \omega_{e})$ is an grading-restricted
conformal open-string vertex algebra satisfying the 
positive weight condition above and with
$V$ in its meromorphic center. Again, the proof is similar to
the proof of the result in \cite{HKL} mentioned above. 
For reader's convenience,
we give a proof  here. 

The identity property for the vacuum follows immediately {}from the left
unit property $\mu\circ (\iota_{V_{e}}\boxtimes \id_{V_{e}})\circ
l^{-1}_{V_{e}} =\id_{V_{e}}$.  The creation property follows {}from the
right unit property 
$\mu \circ (\iota_{V_{e}}\otimes  
\id_{V_{e}})\circ r^{-1}_{V_{e}}=\id_{V_{e}}$.  
The Virasoro relations and the $L(0)$-grading
property follows {}from the fact that $V_{e}$ is a $V$-module. The
$L(-1)$-derivative property and the commutator formula for the Virasoro
operators and $\Y^{f}_{e}$ follow {}from the fact that $\Y^{f}_{e}$ is an
intertwining operator. 

We now prove associativity. As above, for any $r\in \R_{+}$,
let 
$$\mu_{r}: V_{e}\boxtimes_{P(r)}V_{e}\to V_{e}$$
be the module map corresponding to the intertwining operator $\Y^{f}_{e}$. 
By definition, we have
\begin{equation}
\mu_{r}(u\boxtimes_{P(r)} v)=\Y^{f}_{e}(u, r)v=(\mu\circ \mathcal{T}_{\gamma})
(u\boxtimes_{P(r)} v)\label{converse3}
\end{equation}
for $u, v\in V_{e}$ and $r\in \R_{+}$, 
where $\gamma$ is a path 
{}from $r$ to $1$ in $\R_{+}$. 
By definition, for $r_{1}, r_{2}\in \R_{+}$ satisfying
$r_{1}>r_{2}>r_{1}-r_{2}> 0$,  paths 
$\gamma_{1}$ and $\gamma_{2}$ in $\R_{+}$ 
{}from $1$ to $r_{1}$, $r_{2}$, respectively,
and paths $\gamma_{3}$ 
and $\gamma_{4}$ in $\R_{+}$ 
{}from $r_{2}$ and $r_{1}-r_{2}$ to $1$, respectively,
(\ref{mu-assoc3})--(\ref{mu-assoc4}) hold.

Compose both sides of the associativity 
$$\mu \circ(\id \boxtimes \mu)=(\mu \circ(\mu \boxtimes \id))\circ 
\mathcal{A}$$
for the $\mathcal{C}$-algebra $V_{e}$ with 
$$((\id \boxtimes_{P(z_{1})} 
\mathcal{T}_{\gamma_{2}})\circ \mathcal{T}_{\gamma_{1}})^{-1},$$
where $r_{1}, r_{2}\in \R_{+}$ satisfying
$r_{1}>r_{2}>r_{1}-r_{2}>0$ and 
$\gamma_{1}$ and $\gamma_{2}$, as above,
are paths {}from $1$ to $r_{1}$ and $r_{2}$, 
respectively, in $\R_{+}$.
Then we obtain
\begin{eqnarray}
\lefteqn{\mu \circ(\id \boxtimes \mu)\circ
((\id \boxtimes_{P(r_{1})} 
\mathcal{T}_{\gamma_{2}})\circ \mathcal{T}_{\gamma_{1}})^{-1}}\nn
&&=(\mu \circ(\mu \boxtimes \id))\circ 
\mathcal{A}\circ ((\id \boxtimes_{P(r_{1})} 
\mathcal{T}_{\gamma_{2}})\circ \mathcal{T}_{\gamma_{1}})^{-1}.\label{converse4}
\end{eqnarray}
Using (\ref{mu-assoc3})--(\ref{mu-assoc5}) and (\ref{converse4}),
we obtain
\begin{eqnarray}
\lefteqn{\mu_{r_{1}} \circ(\id \boxtimes_{P(r_{1})} \mu_{r_{2}})}\nn
&&=\mu \circ(\id \boxtimes \mu)\circ
((\id \boxtimes_{P(r_{1})} 
\mathcal{T}_{\gamma_{2}})\circ \mathcal{T}_{\gamma_{1}})^{-1}\nn
&&=(\mu \circ(\mu \boxtimes \id))\circ 
\mathcal{A}\circ ((\id \boxtimes_{P(r_{1})} 
\mathcal{T}_{\gamma_{2}})\circ \mathcal{T}_{\gamma_{1}})^{-1}\nn
&&=(\mu_{r_{2}} \circ(\mu_{r_{1}-r_{2}} \boxtimes_{P(r_{2})} \id))\circ
(\mathcal{T}_{\gamma_{3}}\circ (\mathcal{T}_{\gamma_{4}}
\boxtimes_{P(r_{2})} \id))^{-1}\nn
&&\quad \quad\quad\quad\quad\circ
\mathcal{A}\circ ((\id \boxtimes_{P(r_{1})} 
\mathcal{T}_{\gamma_{2}})\circ \mathcal{T}_{\gamma_{1}})^{-1}\nn
&&=(\mu_{r_{2}} \circ(\mu_{r_{1}-r_{2}} \boxtimes_{P(r_{2})} \id))
\circ \mathcal{A}^{P(r_{1}-r_{2}), P(r_{2})}_{P(r_{1}), P(r_{2})}.
\label{converse5}
\end{eqnarray}

For the next step, we use the convergence of products and iterates
of intertwining operators for $V$. Because of the convergence,
$\id \boxtimes_{P(r_{1})} 
\overline{\mu_{r_{2}}}$ is well defined and it is clear that
$\overline{\mu_{r_{1}}} \circ(\id \boxtimes_{P(r_{1})} 
\overline{\mu_{r_{2}}})$ is equal to
$\overline{\mu_{r_{1}} \circ(\id \boxtimes_{P(r_{1})} 
\mu_{r_{2}})}$. Similarly, $\overline{\mu_{r_{1}-r_{2}}}
\boxtimes_{P(r_{2})} \id$ is well-defined and 
$\overline{\mu_{r_{1}}} \circ(\overline{\mu_{r_{1}-r_{2}}}
\boxtimes_{P(r_{2})} \id)$ is equal to
$\overline{\mu_{r_{1}} \circ(\mu_{r_{1}-r_{2}}
\boxtimes_{P(r_{2})} \id)}$.
Taking the natural completions of both sides of 
(\ref{converse5}), we obtain
\begin{equation}
\overline{\mu_{r_{1}}} \circ(\id \boxtimes_{P(r_{1})} 
\overline{\mu_{r_{2}}})
=\overline{\mu_{r_{1}}} \circ(\overline{\mu_{r_{1}-r_{2}}}
\boxtimes_{P(r_{2})} \id)\circ \overline{\mathcal{A}}^{P(r_{1}-r_{2}), 
P(r_{2})}_{P(r_{1}), P(r_{2})}.\label{converse6}
\end{equation}
Applying both sides of (\ref{converse6}) to $u\boxtimes_{P(r_{1})}(v
\boxtimes_{P(r_{2})}w)$ for $u, v, w\in V_{e}$,
pairing the result with $v'\in V_{e}$ and using (\ref{converse3})
and 
$$\overline{\mathcal{A}}^{P(r_{1}-r_{2}), 
P(r_{2})}_{P(r_{1}), P(r_{2})}(u\boxtimes_{P(r_{1})}(v
\boxtimes_{P(r_{2})}w))=(u\boxtimes_{P(r_{1}-r_{2})}v)
\boxtimes_{P(r_{2})}w,$$
we obtain the associativity 
$$\langle v', Y^{O}_{e}(u, r_{1})Y^{O}_{e}(v, r_{2})w\rangle
=\langle v', Y^{O}_{e}(Y^{O}_{e}(u, r_{1}-r_{2})v, r_{2})w\rangle$$
for $u, v, w\in V_{e}$, $v'\in V_{e}'$ and $r_{1}, r_{2}\in \R_{+}$
satisfying $r_{1}>r_{2}>r_{1}-r_{2}>0$.

We now prove that $\iota_{V_{e}}(V)$ is in the meromorphic center of $V_{e}$.
Clearly $\iota_{V_{e}}(V)$ is a vertex operator algebra isomorphic to 
$V$, $\iota_{V_{e}}$ is an isomorphism of vertex operator
algebras {}from $V$ to $\iota_{V_{e}}(V)$ and 
thus $V_{e}$ is an $\iota_{V_{e}}(V)$-module.
We know that the restriction $\Y^{f}_{e}|_{V_{e}\otimes \iota_{V_{e}}(V)}$ of 
$\Y^{f}_{e}$ to $V_{e}\otimes \iota_{V_{e}}(V)$
is in fact the intertwining operator 
of type $\binom{V_{e}}{V_{e}\iota_{V_{e}}(V)}$ for the vertex operator 
algebra $\iota_{V_{e}}(V)$ corresponding to 
the module map $\mu|_{V_{e}\boxtimes \iota_{V_{e}}(V)}: 
V_{e}\boxtimes \iota_{V_{e}}(V)\to V_{e}$
which is the restriction of $\mu$ to $V_{e}\boxtimes \iota_{V_{e}}(V)$.
By the creation property for $Y^{O}_{e}$, we have 
$$\lim_{r\to 0}\Y^{f}_{e}(u, r)\mathbf{1}_{e}
=\lim_{r\to 0}Y^{O}_{e}(u, r)\one_{e}=u$$
for $u\in V_{e}$. Since the space of intertwining operators of type 
$\binom{V_{e}}{V_{e}\iota_{V_{e}}(V)}$ 
is isomorphic to the space of intertwining 
operators of type $\binom{V_{e}}{\iota_{V_{e}}(V)V_{e}}$, which in turn 
is isomorphic to the space of module maps {}from $V_{e}$ to itself,
any intertwining operator $\Y$ of this type satisfying the creation 
property
$$\lim_{r\to 0}\Y(u, r)\mathbf{1}_{e}=u$$ 
must be equal to $\Y^{f}_{e}|_{V_{e}\otimes \iota_{V_{e}}(V)}$. In fact, 
the intertwining operator $\Y$ of such type defined by 
$$\Y(u, x)v=e^{xL(-1)}Y_{V_{e}}(v, -x)u$$
for $u\in V_{e}, v\in \iota_{V_{e}}(V)$, 
where $Y_{V_{e}}$ is the vertex operator map
for the $\iota_{V_{e}}(V)$-module $V_{e}$, is such an intertwining operator.
Thus we have
\begin{equation}\label{mero-cent}
\Y^{f}_{e}|_{V_{e}\otimes \iota_{V_{e}}(V)}(u, x)v
=e^{xL(-1)}Y_{V_{e}}(v, -x)u
\end{equation}
for $u\in V_{e}, v\in \iota_{V_{e}}(V)$. 
But both $Y_{V_{e}}$ and 
$\Y^{f}_{e}|_{\iota_{V_{e}}(V)\otimes V_{e}}$ are intertwining 
operators of type $\binom{V_{e}}{\iota_{V_{e}}(V)V_{e}}$ 
satisfying the identity 
property
and the space of intertwining operators of such type, as we mentioned 
above, is isomorphic
to the space of module maps {}from $V_{e}$ to itself.
So  $Y_{V_{e}}$ 
and $\Y^{f}_{e}|_{\iota_{V_{e}}(V)\otimes V_{e}}$ must be equal. 
Thus (\ref{mero-cent}) says that 
$\iota_{V_{e}}(V)$ is in the meromorphic center of 
$V_{e}$. So $\iota_{V_{e}}$ is an injective homomorphism 
{}from $V$ to the meromorphic center of $V_{e}$. 

The constructions above give two functors and it is easy to see that 
they are inverse to each other. Thus the two categories are
isomorphic.

Finally we prove the last statement. 
We assume that  the intertwining operator algebra on the 
direct sum of all irreducible $V$-modules satisfies the
positive weight condition. In particular, as an 
open-string vertex algebra, $V$ itself satisfies the positive weight 
condition. Let $V_e$ be a grading-restricted conformal 
open-string vertex algebra and $\iota_{V_{e}}$ an injective 
homomorphism of vertex operator algebras {}from $V$ to $C_{0}(V_{e})$. 
Since the weights 
of the nonzero elements of all the 
irreducible $V$-modules are nonnegative, the weights of the
nonzero elements of the $V$-module $V_{e}$ are also nonnegative.
Assume that 
$V_{e}$ satisfies the positive weight condition.
Let $f\in \hom_{\mathcal{C}}(V, V_{e})$. Since $f$ preserves the grading
and since $V$ and $V_{e}$ both satisfy the positive weight condition, 
it is clear that $f$ maps
$\mathbf{1}$ to a scalar multiple of $\mathbf{1}_{e}$. Since 
$V$ as a module is generated by $\mathbf{1}$, $f$ is determined 
completely by the scalar above. On the other hand, given any scalar,
we can also construct an element of $\hom_{\mathcal{C}}(V, V_{e})$
such that it maps $\mathbf{1}$ to the scalar times $\mathbf{1}_{e}$.
Thus $\dim \hom_{\mathcal{C}}(V, V_{e})=1$.  Conversely, 
assume that $\dim \hom_{\mathcal{C}}(V, V_{e})=1$. We already 
know that the weights of 
nonzero elements of the $V$-module $V_{e}$ are also nonnegative.
Assume that there is an element of $(V_{e})_{(0)}$
which is not proportional to 
$\one_{e}$. Then this element  generates a $V$-submodule 
of the $V$-module $V_{e}$. Since all $V$-modules are completely reducible,
we can find an irreducible $V$-submodule of this $V$-submodule
such that it is generated by an element of $(V_{e})_{(0)}$
which is not proportional to 
$\one_{e}$. Since any irreducible $V$-module having a nonzero element of 
weight $0$ must 
be isomorphic to $V$, this $V$-submodule is isomorphic to $V$. But this 
$V$-submodule is not equal to $\iota_{e}(V)\subset C_{0}(V_{e})$ since its 
generator of weight $0$ is not proportional to $\one_{e}$. Thus we see that 
$\dim \hom_{\mathcal{C}}(V, V_{e})>1$. Contradiction. So 
$V_{e}$ satisfies the positive weight condition.
\epf

\begin{rema}
{\rm Recall that a {\it commutative associative algebra  in a braided
tensor category $\mathcal{C}$} or a {\it commutative associative 
$\mathcal{C}$-algebra} is an associative $\mathcal{C}$-algebra satisfying 
$\mu \circ \mathcal{R}=\mu$ ({\it commutativity}), where
$\mathcal{R}$ is the commutativity isomorphism {}from $A\otimes A$
to itself. Let $V$ be a vertex operator algebra as in Theorem \ref{main}
and $\mathcal{C}$ the category of $V$-modules. Then an associative 
$\mathcal{C}$-algebra $V_{e}$ is in general not commutative
In fact, for the 
category $\mathcal{C}$ of modules for $V$, the commutativity isomorphism
$\mathcal{R}$ is characterized by 
\begin{equation}
\overline{\mathcal{R}}(u\boxtimes v)=e^{L(-1)}
\overline{\mathcal{T}}_{\gamma_{+}}
(v\boxtimes_{P(-1)} u)\label{commu-iso}
\end{equation}
where $u, v\in V_{e}$, $\gamma_{+}$ is a path
{}from $-1$ to $1$  in the closed upper half plane without passing through
$0$,
$\mathcal{T}_{\gamma_{+}}$ is the corresponding 
parallel transport isomorphism and $\overline{\mathcal{T}}_{\gamma_{+}}$
is the natural extension of $\mathcal{T}_{\gamma_{+}}$ to the algebraic 
completion $\overline{V_{e}\boxtimes V_{e}}$ of $V_{e}\boxtimes V_{e}$. 
The natural extensions of the left- and right-hand sides of commutativity 
applied to $u\boxtimes v$ for $u, v\in V_{e}$
gives
$\overline{\mu}(\overline{\mathcal{R}}(u\boxtimes v))$ and 
$\overline{\mu}(u\boxtimes v)$, respectively. By the characterization
(\ref{commu-iso})
of $\mathcal{R}$ and the relation between $\mu$ and $\Y_{e}^{f}$, 
the left- and right-hand sides of commutativity are further equal to 
$e^{L(-1)}\Y^{f}_{e}(v, -1)u$ and $Y^{f}_{e}(u, 1)v$, respectively. 
Note that in general $\Y^{f}_{e}(v, -1)u\ne Y^{O}_{e}(v, -1)u$.
So $e^{L(-1)}\Y^{f}_{e}(v, -1)u$ and $Y^{f}_{e}(u, 1)v$ are not equal 
in general. Thus commutativity is not true in general. }
\end{rema}

\renewcommand{\theequation}{\thesection.\arabic{equation}}
\renewcommand{\thethm}{\thesection.\arabic{thm}}
\setcounter{equation}{0}
\setcounter{thm}{0}

\section{A geometric and operadic formulation}

In this section, we give a geometric and operadic formulation
of the notion of grading-restricted 
conformal open-string vertex algebra. For the notion of
open-string vertex algebra and other variations,
we have similar geometric and operadic formulations. 
In the present section, we discuss only
grading-restricted conformal open-string vertex algebras.
We assume that the reader is familiar with the geometric 
and operadic formulation of the notion of vertex operator algebra
given by the first author. See \cite{H1}, \cite{H2}, \cite{H5}, 
\cite{HL1} and \cite{HL2} for details.

We first introduce a geometric partial operad. 
Note that $\hat{\mathbb{H}}$ is analytically 
diffeomorphic to the closed unit disk. We use $\Delta_{a}^{r}$ and
$\bar{\Delta}_{a}^{r}$
to denote  the relatively open upper-half disk
in $\bar{\mathbb{H}}$ and the closed upper-half disk
in $\bar{\mathbb{H}}$, respectively, centered at 
$a\in \R$ with radius $r\in \R_{+}$, that is, 
$\Delta_{a}^{r}=B_{a}^{r}\cap \overline{\HH}$ and
$\bar{\Delta}_{a}^{r}=\bar{B}_{a}^{r}\cap \overline{\HH}$ where 
$B_{a}^{r}$ and $\bar{B}_{a}^{r}$ are the open and closed disks 
centered at $a\in \R$ with radius $r\in \R_{+}$.

A {\it disk with 
strips of type $(m, n)$} ($m, n\in \N$) 
is a disk $S$ (a genus-zero compact connected
one-dimensional complex manifold with one connected component of
boundary)  with $m+n$ distinct, ordered points $p_{1}, \dots , p_{m+n}$
(called {\it boundary punctures}) on the boundary of 
$S$ with $p_{1}, \dots, p_{m}$ negatively
oriented and the other punctures positively oriented, and with local
analytic coordinates 
$$(U_{1}, \varphi _{1}), \dots, (U_{m+n}, \varphi
_{m+n})$$
vanishing at the boundary punctures $p_{1}, \dots , p_{m+n}$,
respectively, where for each $i=1, \dots, m+n$, $U_{i}$ is a local coordinate
neighborhood at $p_{i}$ and $\varphi _{i}: U_{i} \to \bar{\HH}$, mapping
the boundary part of $U_{i}$ analytically to $\R$ and satisfying $\varphi_{i}(p_{i})=0$, 
is a local analytic coordinate map vanishing at $p_{i}$. In the present paper,
we consider only disks with strips of types $(1, n)$ for $n\in \N$. For 
such a disk with strips, we use the subscript $0$ and the subscripts
$1, \dots, n$ to indicate 
that the corresponding boundary punctures are negatively oriented and 
positively oriented, respectively.

Let $S_{1}$ and $S_{2}$ be disks with strips of type $(1, m)$ and of
type $(1, n)$, respectively.  Let $p_{0}, \dots, p_{m}$ be the boundary
punctures of $S_{1}$, $q_{0},\dots, q_{n}$ the boundary punctures of
$S_{2}$, $(U_{i}, \varphi_{i})$ the local coordinate at $p_{i}$ for some
fixed $i$ satisfying $0<i\le m$, and $(V_{0}, \psi_{0})$ the local
coordinate at $q_{0}$. Note that in our convention discussed above,
$p_{0}$ and $q_{0}$ are the negatively oriented boundary punctures on
$S_{1}$ and $S_{2}$, respectively. Assume that there exists $r\in
\R_{+}$ such that $\varphi _{i}(U_{i})$ contains $\bar{\Delta}_{0}^{r}$
and $\psi_{0}(V_{0})$ contains $\bar{\Delta}_{0}^{1/r}$.  Assume also
that $p_{i}$ and $q_{0}$ are the only boundary punctures in
$\varphi_{i}^{-1}(\bar{\Delta}_{0}^{r})$ and
$\psi_{0}^{-1}(\bar{\Delta}_{0}^{1/r})$, respectively. In this case we
say that {\it the $i$-th boundary puncture of the first disk with strips
can be sewn with the $0$-th boundary puncture of the second disk with
strips}. {}From these two disks with strips we obtain a disk with strips
of type $(1, m+n-1)$ by cutting $\varphi_{i}^{-1}(\Delta_{0}^{r})$ and
$\psi_{0}^{-1}(\Delta_{0}^{1/r})$ {}from $S_{1}$ and $S_{2}$,
respectively, and then identifying the new parts of the boundaries (the
parts not on the boundaries of the original surfaces) of the resulting
surfaces using the map $\varphi_{i}^{-1} \circ (-J) \circ \psi_{0}$
where $J$ is the map {}from $\mathbb{C}^{\times}$ to itself given by
$J(w)=1/w$.  The boundary punctures (with ordering) of this disk with
strips are $p_{0}, \dots,p_{i-1}$, $q_{1}, \dots, q_{n}$, $p_{i+1},
\dots,p_{m}$. The local coordinates vanishing at these punctures are
given in the obvious way. This sewing procedure gives a partial
operation which we call the {\it sewing operation}. Note that we have to
use $-J$ instead of $J$ (as in \cite{H5}) in the definition of the
sewing operation. 

We define the notion of conformal equivalence between two disks with
strips in the obvious way. The space of equivalence classes of disks
with strips is called the {\it moduli space of disks with strips}.  
Similar to  the moduli spaces of spheres with tubes in 
\cite{H5}, the
moduli space of disks with strips of type $(1, n)$ ($n\ge 1$) can be
identified with $\Upsilon(n)=\Lambda^{n-1}\times \Pi \times
\Pi_{\R_{+}}^{n}$ where $\Pi$ is the set of all sequences 
$A=\{A_{j}\}_{j\in \Z_{+}}$ of
real numbers such that $$\mbox{\rm exp}\left(\displaystyle
\sum_{j>0}A_{j}x^{j+1}\frac{d}{dx}\right) x$$ is a convergent power
series in some neighborhood of $0$, $\Pi_{\R_{+}}=
\R_{+} \times \Pi$, and $\Lambda^{n-1}$ is the set of
elements of $\R^{n-1}$ with nonzero and distinct components. We
think of each element of $\Upsilon(n)$, $n\ge 1$, as the disk $\hat{\HH}$
equipped with ordered punctures $\infty$, $r_{1},\dots,r_{n-1}$, $0$,
with an element of $\Pi$ specifying the local coordinate at $\infty$ and
with $n$ elements of $\Pi_{\R_{+}}$ specifying the local coordinates
at the other punctures. Analogously, the moduli space of disks with
strips of type $(1, 0)$ can be identified with $\Upsilon(0)=\{A\in
\Pi\;|\;A_{1}=0\}$. Then the moduli space of disks with strips can be
identified with $\cup_{n\ge 0}\Upsilon(n)$.  {}From now on we will refer to
$\cup_{n\in \N}\Upsilon(n)$ as the moduli space of disks with strips.  The
sewing operation for disks with strips induces a partial operation on
$\cup_{n\in \N}\Upsilon(n)$.  It is still called the sewing operation.

Let $I_{\Upsilon}\in \Upsilon(1)$ be the equivalence class 
containing the standard disk
$\hat{\HH}$ with  the negatively oriented puncture $\infty$, the only
positively oriented puncture $0$ , and with standard local coordinates
vanishing at $\infty$ and $0$. Here for $a\in \R\subset \hat{\HH}$, 
the standard local coordinate vanishing at $a$ is given by 
$w\mapsto  w-a$.
and for $\infty\in \hat{\HH}$, the standard local coordinate 
vanishing at $\infty$ is given by $w\mapsto -\frac{1}{w}$. 
Note the minus sign in the definition of the standard local 
coordinate at $\infty$.
For $n\in \mathbb{N}$, 
the symmetric group $S_{n}$ acts on $\Upsilon(n)$ in an obvious way.
Then by construction, the following result is clear:

\begin{prop}
The sequences $\Upsilon=\{\Upsilon(n)\;|\;n\in \mathbb{N}\}$ of 
moduli spaces,
together with the sewing operation, the identity $I_{\Upsilon}$ and the 
actions of the symmetric groups, has a structure of
an associative smooth $\R_{+}$-rescalable partial operad.
\end{prop}

We shall call the $\R_{+}$-rescalable 
partial operad $\Upsilon$ the {\it boundary disk partial operad}. 
Note that the boundary disk partial operad is very different {}from the 
so-called little disk operad which are constructed using 
the embeddings of disks in the unit disk. 
In fact, $\Upsilon$ can be viewed as a partial suboperad of 
the sphere partial operad $K$ discussed in \cite{H5}. 
Geometrically, any disk with strips of type $(1, n)$ is conformally 
equivalent to a disk with strips of type $(1, n)$
whose underlying disk is $\hat{\HH}$ and whose negatively oriented 
puncture is $\infty$. But any such disk with strips of type $(1, n)$
corresponds to a sphere with tubes of type $(1, n)$ whose 
underlying sphere is $\hat{\mathbb{C}}$, whose punctures
are the same as  those on the disk with strips, whose local 
coordinates
vanishing at positively oriented punctures are
the analytic extensions
of those on the disk with strips and whose local coordinate 
vanishing at the negatively oriented puncture is the analytic extension
of the negation of that on the disk with strips. 
Thus we obtain a map {}from 
$\Upsilon(n)$ to $K(n)$ and this map is clearly injective. 
In fact the images of $\Upsilon(n)$ in $K(n)$ for $n\ge 2$ are
\begin{eqnarray*}
\lefteqn{\{(r_{1}, \dots, r_{n-1}; A^{(0)}, (a_{0}^{(1)}, A^{(1)}),
\cdots, (a_{0}^{(n)}, A^{(n)}))\in K(n)\;|}\nn
&&\quad\quad \quad\quad\quad\quad r_{1}, \dots, r_{n-1}\in \R, 
a_{0}^{(1)}, \dots, a_{0}^{(n)}\in 
\R_{+}, A^{(0)}, \dots, A^{(n)}\in \Pi\}.
\end{eqnarray*}
The images of $\Upsilon(0)$ in $K(0)$ and of $\Upsilon(1)$ in $K(1)$ are 
$$\{A^{(0)}\in K(0)\;|\;  A^{(0)}\in \Pi,\; A^{(0)}_{1}=0\}$$
and 
$$\{(A^{(0)}, (a_{0}^{(1)}, A^{(1)}))\in K(1)\;|\; 
a_{0}^{(1)}\in 
\R_{+}, A^{(0)}, A^{(1)}\in \Pi\},$$
respectively. 
In addition, by the definitions of the maps {}from $\Upsilon(n)$ to 
$K(n)$ for $n\in \N$ and the sewing operations in $\Upsilon$ and $K$, 
it is clear that the maps {}from $\Upsilon(n)$
to $K(n)$ for $n\in \mathbb{N}$ respect 
the sewing operations, the identities and 
the actions of $S_{n}$ and thus give an injective morphism of 
partial operads.
{}From now on, we shall identify the partial operad $\Upsilon$ with 
its image in $K$ under this injective morphism. 

For any $c\in \mathbb{C}$, the restriction of the 
partial operad $\tilde{K}^{c}$ of the $\frac{c}{2}$-th power of 
the determinant line bundles over $K$ to 
$\Upsilon$ gives a partial suboperad $\tilde{\Upsilon}^{c}$ of 
$\tilde{K}^{c}$. This partial operad is called the {\it 
$\mathbb{C}$-extension
of $\Upsilon$ of central charge $c$}. 

We now consider certain (pseudo-)algebras over the partial operad 
$\tilde{\Upsilon}^{c}$ for $c\in \C$. In the terminology of \cite{HL1}, 
\cite{HL2} and \cite{H5}, we consider $\tilde{\Upsilon}^{c}$-associative 
(pseudo-)algebras
satisfying an additional differentiability condition. Since the rescaling 
group of $\tilde{\Upsilon}^{c}$ is $\R_{+}$, we need to consider 
modules for $\R_{+}$.
Since an equivalence class of irreducible modules for 
$\R_{+}$ is determined by a real number $s$ such that $a\in 
\R_{+}$ acts on modules in this class as the scalar multiplication by
$a^{-s}$, any completely reducible module for $\R_{+}$ is of
the form $V=\coprod_{s\in \R}V_{(s)}$ where $V_{(s)}$ is the sum
of the $\R_{+}$-submodules in the class determined by the
real number $s$.  We shall consider only those algebras 
over $\tilde{\Upsilon}^{c}$ whose underlying vector space  is of the form 
$V=\coprod_{s\in \R}V_{(s)}$
such that $\dim V_{(s)}<\infty$. Recall {}from \cite{HL1}, \cite{HL2}
and \cite{H5} that given any $\R_{+}$-submodule
$W$ of $V$, the endomorphism 
partial pseudo-operad associated to the pair $(V, W)$ is the sequence
$H^{\R_{+}}_{V, W}=\{H^{\R_{+}}_{V, W}(n)\}_{n\in \N}$.
where $H^{\R_{+}}_{V, W}(n)$ is the set of all multilinear maps {}from 
$V^{\otimes n}$ to $\overline{V}$ such that $W^{\otimes n}$ is mapped to 
$\overline{W}$, equipped with natural operadic structures.

\begin{defn}
{\rm A   {\it differentiable (or $C^{1}$)
$\tilde{\Upsilon}^{c}$-associative pseudo-algebra} 
is a completely reducible $\R_{+}$-module $V=\coprod_{s\in \R}V_{(s)}$
satisfying the condition $\dim V_{(r)}<\infty$ for $s\in \R$ 
equipped with an $\R_{+}$-submodule
$W$ and a morphism $\Phi$ of $\R_{+}$-rescalable pseudo-partial operad
{}from $\tilde{\Upsilon}^{c}$ to the endomorphism 
partial pseudo-operad $H^{\R_{+}}_{V, W}$ (that is, an 
$\tilde{\Upsilon}^{c}$-associative pseudo-algebra) 
satisfying the following conditions:

\begin{enumerate}

\item For $s$ sufficiently negative, $V_{(s)}=0$. 

\item For any $n\in \N$, $\Phi_{n}: \tilde{\Upsilon}^{c}(n)\to 
H^{\R_{+}}_{V, W}(n)$ is linear on the fibers of 
$\tilde{\Upsilon}^{c}(n)$.

\item For any $s_{1}, \dots, s_{n}\in \R$, there exists a finite 
subset $R(s_{1}, \dots, s_{n})\subset \R$ such that 
the image of $\coprod_{s\in s_{1}+\Z}V_{(s)}\otimes \cdots
\otimes \coprod_{s\in s_{n}+\Z}V_{(s)}$ under 
$\Phi_{n}(\psi_{n}(Q))$ for any $Q\in \tilde{\Upsilon}^{c}(n)$
is in $\coprod_{s\in R(s_{1}, \dots, s_{n})+\Z}V_{(s)}$.

\item For any $v'\in V'$, $v_{1}, \dots, v_{n}\in V$, $\langle v',
\Phi_{n}(\psi_{n}(Q))(v_{1}\otimes \dots\otimes
v_{n})\rangle$  as a function of
$$Q=(r_{1}, \dots, r_{n-1}; A^{(0)}, (a_{0}^{(1)}, A^{(1)}),
\cdots, (a_{0}^{(n)}, A^{(n)}))\in \tilde{\Upsilon}^{c}(n)$$
is of the 
form
$$\sum_{i=1}^{m}f_{i}(r_{1}, \dots, r_{n-1})
g_{i}(A^{(0)}, (a_{0}^{(1)}, A^{(1)}),
\cdots, (a_{0}^{(n)}, A^{(n)}))$$
where $f_{i}(r_{1}, \dots, r_{n-1})$ for $i=1, \dots, m$ 
are continuous differentiable functions 
of $r_{1}, \dots, r_{n-1}$ and $g_{i}(A^{(0)}, (a_{0}^{(1)}, A^{(1)}),
\cdots, (a_{0}^{(n)}, A^{(n)}))$ for 
$i=1, \dots, m$ are polynomials in $A^{(0)}, (a_{0}^{(1)})^{\pm 1}, 
A^{(1)},
\cdots, (a_{0}^{(n)})^{\pm 1}, A^{(n)}$. 

\end{enumerate}

{\it Morphisms} (respectively, {\it isomorphisms}) of 
differentiable $\tilde{\Upsilon}^{c}$-associative pseudo-algebras
are morphisms (respectively, isomorphisms) of the
underlying $\tilde{\Upsilon}^{c}$-associative pseudo-algebras. }
\end{defn}

We denote the differentiable $\tilde{\Upsilon}^{c}$-associative
pseudo-algebra just defined by $(V, W, \Phi)$ or simply $V$. It is easy
to see that a differentiable $\tilde{\Upsilon}^{c}$-associative
pseudo-algebra is actually analytic in the sense that for any $v'\in
V'$, $v_{1}, \dots, v_{n}$ $\in V$, $\langle v', \nu(Q)(v_{1}, \dots,
v_{n})\rangle$ is analytic in $Q$ because of the sewing axiom (that is,
the sewing operation in $\Upsilon$ corresponds to the contraction in
$H^{\R_{+}}_{V, W}$ under $\Phi$). Using this fact and the fact that the
expansion of analytic functions are always absolutely convergent in the
domain of convergence, it is easy to obtain:

\begin{prop}\label{pseudo}
Any differentiable $\tilde{\Upsilon}^{c}$-associative pseudo-algebra 
$(V, W, \Phi)$ is an 
$\tilde{\Upsilon}^{c}$-associa\-tive algebra, that is, the image 
of $\tilde{\Upsilon}^{c}$ under $\Phi$ is a partial operad
(the image 
of $\tilde{\Upsilon}^{c}$ under $\Phi$
satisfies  the composition-associativity).
\end{prop}

We omit the proof of this result since it is the same as 
the proof of the corresponding result in \cite{H5}. Because 
of this result, we shall call a differentiable 
$\tilde{\Upsilon}^{c}$-associative pseudo-algebra simply
a {\it differentiable 
$\tilde{\Upsilon}^{c}$-associative algebra}.

Now we have the following main theorem which gives a geometric 
and operadic formulation of the notion of grading-restricted
conformal open-string vertex algebras:

\begin{thm}\label{main1}
The category of grading-restricted conformal open-string vertex algebras
of central charge
$c$ is isomorphic to the category of  
differentiable $\tilde{\Upsilon}^{c}$-associative algebras.
\end{thm}
\pf 
The proof of this theorem is basically the same as that of the 
isomorphism theorem for the geometric and operadic formulation of 
vertex operator algebras in \cite{H5}. Here we give a sketch.
Some more details will be given in \cite{K}. 

Let $(V, Y^{O}, \mathbf{1}, \omega)$ be 
a grading-restricted conformal open-string vertex algebra
of central charge $c$. We construct a
differentiable $\tilde{\Upsilon}^{c}$-associative algebras
of central charge $c$ as follows:  The $\R$-graded
vector space $V$ is naturally a completely reducible 
$\R_{+}$-module.  The module $W$ for the Virasoro algebra
generated by $\one$ is an $\R$-graded subspace of $V$ and
therefore is an $\R_{+}$-submodule of $V$. In \cite{H2} and
\cite{H5}, a  section $\psi$ of the line bundle $\tilde{K}^{c}$ over
$K$ is chosen.  The restriction of this section to $\Upsilon$ is a
section of $\tilde{\Upsilon}^{c}$ and, for simplicity, 
we still denote it by 
$\psi$. For an element
\begin{equation}\label{Q}
Q=(r_{1}, \dots, r_{n-1}; A^{(0)}, (a^{(1)}_{0}, A^{(1)}), 
\dots, (a^{(n)}_{0}, A^{(n)}))
\end{equation}
of $\Upsilon(n)$,
any element of the fiber of $\tilde{\Upsilon}^{c}$
over $Q$ is of the
form $\lambda\psi_{n}(Q)$ where $\lambda \in \mathbb{C}$. When 
$r_{1}>\cdots >r_{n-1}>0$, we define
$\Phi_{n}(\lambda \psi_{n}(Q))$ by
\begin{eqnarray*}
\lefteqn{(\Phi_{n}(\lambda \psi_{n}(Q)))(v_{1}\otimes \cdots \otimes 
v_{n})}\nn
&&=\lambda e^{-\sum_{j\in \Z_{+}}A^{(0)}_{j}L(-j)} 
Y^{O}(e^{-\sum_{j\in \Z_{+}}A^{(1)}_{j}L(j)}(a^{(1)}_{0})^{-L(0)}v_{1}, 
r_{1})\cdots\nn
&&\hspace{6em}\cdot
Y^{O}(e^{-\sum_{j\in \Z_{+}}A^{(n-1)}_{j}L(j)}(a^{(n-1)}_{0})^{-L(0)}v_{n-1}, 
r_{n-1}) 
\cdot\nno\\
&&\hspace{6em}\cdot 
e^{-\sum_{j\in \Z_{+}}A^{(n)}_{j}L(j)}(a^{(n)}_{0})^{-L(0)}
v_{n}
\end{eqnarray*} 
for $v_{1}, \dots, v_{n}\in V$. In general, for any $Q\in \Upsilon(n)$,
we can always find 
$\sigma_{Q}\in S_{n}$ such that $\sigma_{Q}(Q)$ is of the form 
of the right-hand side of (\ref{Q}) such 
that $r_{1}>\cdots >r_{n-1}>0$. We define 
$\Phi_{n}(\lambda \psi_{n}(Q))$ by
$$(\Phi_{n}(\lambda \psi_{n}(Q)))(v_{1}\otimes \cdots \otimes 
v_{n})=\Phi_{n}(\lambda\psi_{n}(\sigma_{Q}(Q)))(v_{\sigma_{Q}^{-1}(1)}\otimes
\cdots \otimes v_{\sigma_{Q}^{-1}(n)})$$
for $v_{1}, \dots, v_{n}\in V$.
It can be  verified in the same way
as in \cite{H5}  that the triple $(V,W,\nu)$ is 
a differentiable $\tilde{\Upsilon}^{c}$-associative algebra
of central charge $c$. This construction
gives a functor {}from the category of grading-restricted 
conformal open-string vertex algebras
of central charge $c$ to the category of 
differentiable $\tilde{\Upsilon}^{c}$-associative algebras.

Conversely, given a differentiable $\tilde{\Upsilon}^{c}$-associative 
algebra $(V,W,\Phi)$, we construct a grading-restricted 
conformal open-string vertex algebra 
as follows: As in \cite{H5}, for $\varepsilon\in \R$ and $i\in \Z_{+}$, 
let $A(\varepsilon; i)$ be 
the element of $\Pi$ whose 
$i$-th component is equal to $\varepsilon$ and all other components
are $0$ and $\mathbf{0}$ the element of $\Pi$ whose components
are all $0$, and for $r\in \R_{+}$,
let 
$$P(r)=(r; \mathbf{0}, (1, \mathbf{0}), (1, \mathbf{0}))
\in \Upsilon(2)\subset K(2).$$
We define the vertex operator map
\begin{eqnarray*}
Y^{O}: (V\otimes V)\times \mathbb{R}_{+}&\to& \overline{V},\nn
(v_{1}\otimes v_{2}, r)&\mapsto& Y^{O}(v_{1}, r)v_{2}
\end{eqnarray*}
by 
$$Y^{O}(v_{1}, r)v_{2}=(\Phi_{2}(\psi_{2}((P(r))))(v_{1}\otimes v_{2})$$
for $v_{1}, v_{2}\in V$ and $r\in \R_{+}$.
The vacuum $\mathbf{1}\in V$ is given by 
$$\mathbf{1}=\Phi_{0}(\psi_{0}(\mathbf{0})).$$
The conformal element $\omega_{\nu}$ is given by 
$$\omega=-{\displaystyle
\frac{d}{d\varepsilon}}\Phi_{0}(\psi_{0}((A(\varepsilon; 2))))
\lbar_{\varepsilon =0}.$$
It can be proved in the same way
as in \cite{H5} that $(V, Y^{O}, \mathbf{1}, 
\omega)$ is a grading-restricted conformal open-string vertex algebra. 
This construction
gives a functor {}from the category of 
differentiable $\tilde{\Upsilon}^{c}$-associative (pseudo-)algebras
to the category of conformal open-string vertex algebras
of central charge $c$.

It can be  shown in the same way
as in \cite{H5} that these two functors constructed above
are inverse to each other. Thus the conclusion of the theorem is true.
\epfv

The result above can actually be generalized to show 
that a grading-restricted conformal open-string vertex algebra
of central charge $c$
gives an algebra over a partial operad extending the operad
of the $c$-th power of the determinant line bundles over the so-called 
``Swiss-cheese'' operad (see \cite{V}). We actually have a stronger 
isomorphism theorem than Theorem \ref{main1} involving 
meromorphic centers of 
grading-restricted conformal open-string vertex algebras.
To formulate this 
result, we first introduce the underlying  partial operads.

A {\it disk with 
strips and tubes of type $(m, n; k, l)$} ($m, n, k, l\in \N$) 
is a disk $S$ with $m+n$ distinct, ordered points 
$p^{B}_{1}, \dots , p^{B}_{m+n}$
(called {\it boundary punctures}) on the boundary of $S$ 
and $k+l$ distinct, ordered points $p^{I}_{1}, \dots , p^{I}_{k+l}$ 
(called {\it interior punctures}) in the interior of $S$
with $p^{B}_{1}, \dots, p^{B}_{m}$ and $p^{I}_{1}, \dots, p^{I}_{k}$ 
negatively
oriented and the other (boundary or interior) 
punctures positively oriented, and with local
analytic coordinates 
$$(U^{B}_{1}, \varphi^{B}_{1}), \dots, (U^{B}_{m+n}, \varphi^{B}_{m+n}), 
(U^{I}_{1}, \varphi^{I}_{1}), \dots, (U^{I}_{k+l}, \varphi^{i}_{k+l})$$
vanishing at the (boundary or interior) 
punctures $p^{B}_{1}, \dots , p^{B}_{m+n}$, 
$p^{I}_{1}, \dots , p^{I}_{k+l}$,
respectively, where for each $i=1, \dots, m+n$ (or $j=1, \dots, k+l$),
$U^{B}_{i}$ (or $U^{I}_{j}$) is a local coordinate
neighborhood at $p^{B}_{i}$ (or $p^{I}_{j}$)
and $\varphi^{B}_{i}: U^{B}_{i} \to \bar{\HH}$ 
(or $\varphi^{I}_{j}: U^{I}_{j} \to \C$), mapping
the boundary part of $U^{B}_{i}$ (or mapping $U^{I}_{j}$)
analytically to $\R$ (or $\C$) and satisfying $\varphi^{B}_{i}(p^{B}_{i})=0$
(or $\varphi^{I}_{i}(p^{I}_{i})=0$), 
is a local analytic coordinate map vanishing at $p^{B}_{i}$
(or $p^{I}_{i}$). Note that when $k=l=0$, we have 
a disk with strips of type $(m, n)$. In the present paper,
we consider only disks with strips and tubes 
of types $(1, n; 0, l)$ for $n, l\in \N$. For 
such a disk with strips, we use the subscript $0$ and the subscripts
$1, \dots, n$ to indicate 
that the corresponding boundary punctures are negatively oriented and 
positively oriented, respectively.

Similar to disks with strips, we have a sewing operation which
sews two disks with strips and tubes at boundary punctures
of opposite orientations. Here we shall call this sewing operation
the {\it boundary sewing operation}. On the other hand, we can also sew the 
negatively oriented puncture of a sphere with tubes to 
an interior puncture of a disk with strips and tubes
just as we sew two spheres with tubes in \cite{H5}. We shall call 
this sewing operation the {\it interior sewing operation}.

The conformal equivalences for these disks with strips and tubes 
are defined in the obvious way. For $n\ge 1$ and 
$l\in \N$, the moduli space of 
disks with strips and tubes of type $(1, n; 0, l)$ can be identified 
with $\Upsilon(n; l)=\Lambda^{n-1}\times \Pi\times 
\Pi_{\R_{+}}^{n}\times M_{\HH}^{l}
\times H_{c}^{l}$ where $\Lambda$, $\Pi$ and $\Pi_{\R_{+}}$ 
are defined above, $H$ ands $H_{c}$ are defined 
in \cite{H5} and $M_{\HH}^{l}$ is the set of
elements of $\HH^{l}$ with  nonzero and distinct components. 
Analogously, for $l\in \N$, the moduli space of disks with
strips and tubes of type $(1, 0; 0, l)$ can be identified with 
$\Upsilon(0; 1)=\{A\in
\Pi\;|\;A_{1}=0\} \times H_{c}^{l}$. Note that 
$\Upsilon(n)=\Upsilon(n; 0)$ for $n\in \N$. In particular, 
the identity $I_{\Upsilon}$ 
is an element of $\Upsilon(1; 0)$. Also, for $n\in \N$,
$S_{n}$ acts on $\Upsilon(n; l)$ in the obvious way.
Let $\mathfrak{S}(n)=\cup_{l\in \N}\Upsilon(n; l)$ for $n\in \N$.
Then $S_{n}$ acts on $\mathfrak{S}(n)$ for $n\in \N$.
The following result is clear:

\begin{prop}
The sequences $\mathfrak{S}=\{\mathfrak{S}(n)\}_{n\in \N}$
together with the boundary sewing operation, 
the identity $I_{\Upsilon}\in \Upsilon(1)
=\Upsilon(1; 0)$ and the 
actions of the symmetric groups, has a structure of
a smooth $\R_{+}$-rescalable partial operad.
In addition, for each $n\in \N$, there 
is an action of the sphere partial operad $K$ on $\mathfrak{S}(n)$
given by the interior sewing operation. 
\end{prop}

Borrowing the terminology used by Voronov in \cite{V}, 
we shall call the $\R_{+}$-rescalable 
partial operad $\mathfrak{S}$ the {\it Swiss-cheese partial operad}. 
But note that our partial operad is much more larger than the 
Swiss cheese operad. In fact, Swiss cheese operad is an analogue
of the little disk operad while our Swiss cheese partial operad 
is an analogue of the sphere partial operad in \cite{H5}. 

For each pair $n, l\in \N$, 
we have an injective map {}from $\Upsilon(n; l)$ to 
$K(n+2l)$ obtained by doubling disks with strips and tubes
as follows: For any disk with strips and tubes 
of type $(1, n; 0, l)$, by the uniformization theorem, 
we can find a conformally equivalent
disk with strips and tubes of the 
same type such that its underlying disk is $\hat{\HH}$. 
This latter disk with strips and tubes can be doubled to obtain
a sphere with tubes of type $(1, n+2l)$ such that its underlying sphere
is the double $\C\cup \{\infty\}$ of $\hat{\HH}$. 
By definition, we see that
conformally equivalent disks with strips give conformally equivalent
spheres with tubes. Thus we obtain a map {}from $\Upsilon(n; l)$ to
$K(n+2l)$. Clearly this map is injective. 

It is  clear {}from the definition that
these maps respect the (boundary) sewing operations.
In addition, these maps also intertwine the actions of 
$K$ on $\mathfrak{S}(n)$ for 
$n\in \N$ and the actions of $K$ 
on the images of $\mathfrak{S}(n)$ obtained by doubling
the actions of $K$ on $K(n)$.  We shall identify 
$\Upsilon(n; l)$ with its image in 
$K(n+2l)$. 

For any $c\in \mathbb{C}$, the restriction  of the 
partial operad $\tilde{K}^{c}$ of the $\frac{c}{2}$-th power of 
the determinant line bundles over $K$ to 
$\mathfrak{S}$ has a natural structure of a partial  operad. 
This partial operad is called the {\it $\mathbb{C}$-extension
of $\mathfrak{S}$ of central charge $c$}
and is denoted $\tilde{\mathfrak{S}}^{c}$. 
For any $n\in \N$, the action of $K$ on 
$\mathfrak{S}(n)$ also induces an action of 
$\tilde{K}^{c}$ on $\tilde{\mathfrak{S}}^{c}(n)$. For any $n\in \N$, 
the restrictions of the sections 
$\psi_{n+2l}$ of $\tilde{K}^{c}(n+2l)$ for $l\in \N$ to $\Upsilon(n; l)$
gives a section of $\tilde{\mathfrak{S}}^{c}(n)$ and 
we shall use $\psi_{n}^{\mathfrak{S}}$
to denote this section.

We now consider a completely reducible $\R_{+}$-module 
or, equivalently, an $\R$-graded vector spaces $V^{O}=\coprod_{s\in
\R}V^{O}_{(s)}$ and  completely reducible $\C^{\times}$-modules 
or, equivalently, $\Z$-graded vector spaces $V^{LC}=\coprod_{m\in
\mathbb{Z}}V^{LC}_{(m)}$ and $V^{RC}=\coprod_{m\in
\mathbb{Z}}V^{RC}_{(m)}$. (Here $O$, $LC$ and $RC$ means open, 
left closed and right closed, respectively.)
Let $W^{O}$, $W^{LC}$ and $W^{RC}$ be
an $\R_{+}$-submodule
of $V^{O}$, a $\C^{\times}$-submodule of $V^{LC}$ and 
a $\C^{\times}$-submodule of $V^{RC}$, respectively.
Associated to $V^{LC}$, $W^{LC}$, $V^{RC}$, $W^{RC}$, we have the 
endomorphism partial pseudo-operads $H^{\C^{\times}}_{V^{LC}, W^{LC}}$,
$H^{\C^{\times}}_{V^{RC}, W^{RC}}$ and 
$H^{\C^{\times}}_{V^{LC}\otimes (V^{RC})^{-}, W^{LC}\otimes 
(W^{RC})^{-}}$
(see
\cite{HL1}, \cite{HL2} and \cite{H5}), where $(V^{RC})^{-}$ and 
$(W^{RC})^{-}$ are the complex conjugate of $V^{RC}$ and 
$W^{RC}$. We also need an
endomorphism partial operad constructed {}from $V^{O}$, $W^{O}$,
$V^{LC}$, $W^{LC}$, $V^{RC}$ and $W^{RC}$. For $n, l\in \N$, 
let $H^{\R_{+}}_{V^{O}, W^{O}; V^{LC} \otimes (V^{RC})^{-}, W^{LC}\otimes 
(W^{RC})^{-}}(n; l)$
be the space of all linear maps {}from $(V^{O})^{\otimes n}\otimes 
(V^{LC}\otimes (V^{RC})^{-})^{\otimes l}$ to $\overline{V^{O}}$ such that 
$(W^{O})^{\otimes n}\otimes 
(W^{LC}\otimes (W^{RC})^{-})^{\otimes l}$ is mapped to $\overline{W^{O}}$ 
and for $n\in \N$,
let 
\begin{eqnarray*}
\lefteqn{H^{\R_{+}}_{V^{O}, W^{O}; V^{LC} \otimes (V^{RC})^{-}, W^{LC}\otimes 
(W^{RC})^{-}}(n)}\nn
&&=\coprod_{l\in \N}
H^{\R_{+}}_{V^{O}, W^{O}; V^{LC} \otimes (V^{RC})^{-}, W^{LC}\otimes 
(W^{RC})^{-}}(n; l).
\end{eqnarray*}
Then it is clear that for $n\in \N$, the endomorphism partial 
pseudo-operad
$H^{\C^{\times}}_{V^{LC}\otimes 
(V^{RC})^{-}, W^{LC}\otimes (W^{RC})^{-}}$
acts on $H^{\R_{+}}_{V^{O}, W^{O};  V^{LC} \otimes (V^{RC})^{-}, 
W^{LC}\otimes (W^{RC})^{-}}(n)$  and 
\begin{eqnarray*}
\lefteqn{H^{\R_{+}}_{V^{O}, W^{O}; V^{LC} \otimes (V^{RC})^{-}, W^{LC}\otimes 
(W^{RC})^{-}}}\nn
&&=\{H^{\R_{+}}_{V^{O}, W^{O}; V^{LC} \otimes (V^{RC})^{-}, W^{LC}\otimes 
(W^{RC})^{-}}(n)\}_{n\in \N}
\end{eqnarray*}
is an $\R_{+}$-rescalable 
partial pseudo-operad. We call it
the {\it endomorphism partial pseudo-operad for}
$(V^{O}, W^{O}; V^{LC} \otimes (V^{RC})^{-}, W^{LC} \otimes 
(W^{RC})^{-})$.

Notice that 
$H^{\C^{\times}}_{V^{LC}, W^{LC}}\otimes 
(H^{\C^{\times}}_{V^{RC},  W^{RC}})^{-}$ 
(here $(H^{\C^{\times}}_{V^{RC},  W^{RC}})^{-}$ is 
the complex conjugate of $H^{\C^{\times}}_{V^{RC},  W^{RC}}$)
can be embedded naturally into 
the space $H^{\C^{\times}}_{V^{LC}\otimes 
(V^{RC})^{-}, W^{LC}\otimes (W^{RC})^{-}}$. Below we shall view 
$H^{\C^{\times}}_{V^{LC}, W^{LC}}\otimes 
(H^{\C^{\times}}_{V^{RC},  W^{RC}})^{-}$  as a partial pseudo-suboperad of 
$H^{\C^{\times}}_{V^{LC}\otimes 
(V^{RC})^{-}, W^{LC}\otimes (W^{RC})^{-}}$.

Let $\bar{c}$ be the complex conjugate of $c\in \C$. 
The complex conjugate $\overline{\tilde{K}^{\bar{c}}}$ of 
$\tilde{K}^{\bar{c}}$ is also a $\C^{\times}$-rescalable
partial operad. Consequently the tensor product
$\tilde{K}^{c}\otimes \overline{\tilde{K}^{\bar{c}}}$ (the tensor product of
line bundles) is also a $\C^{\times}$-rescalable partial operad. 
Interpreting the action of $K$ on $\mathfrak{S}$
using the method of doubling disks, we see that 
$\tilde{K}^{c}\otimes \overline{\tilde{K}^{\bar{c}}}$ acts naturally on 
$\tilde{\mathfrak{S}}^{c}$. 

We are interested in certain 
algebras over $\tilde{\mathfrak{S}}^{c}$ for $c\in \C$. 

\begin{defn}\label{sw-ch-alg}
{\rm A   {\it  pseudo-algebra
over $\tilde{\mathfrak{S}}^{c}$ generated by a differentiable 
$\tilde{\Upsilon}^{c}$-associative pseudo-algebra
and meromorphic actions of two $\tilde{K}^{c}$-associative algebras} 
or simply a {\it differentiable-meromorphic pseudo-algebra
over $\tilde{\mathfrak{S}}^{c}$}
consists of the following data:

\begin{enumerate}

\item A completely reducible $\R_{+}$-module 
$V^{O}=\coprod_{s\in \R}V^{O}_{(s)}$
satisfying the condition $\dim V^{O}_{(s)}<\infty$ for $s\in \R$
and completely reducible $\C^{\times}$-modules 
$V^{LC}=\coprod_{m\in \mathbb{Z}}V^{LC}_{(m)}$ and 
$V^{RC}=\coprod_{m\in \mathbb{Z}}V^{RC}_{(m)}$
satisfying the condition $\dim V^{LC}_{(m)}<\infty$ and 
$\dim V^{RC}_{(m)}<\infty$ for $m\in \Z$.

\item An $\R_{+}$-submodule
$W^{O}$ of $V^{O}$ and $\C^{\times}$-submodules $W^{LC}$ and 
$W^{RC}$ of $V^{LC}$ and $V^{RC}$, respectively.

\item A morphism $\Phi$ of $\R_{+}$-rescalable partial pseudo-operads {}from 
$\tilde{\mathfrak{S}}^{c}$ to the endomorphism 
partial pseudo-operad $H^{\R_{+}}_{V^{O}, W^{O}; 
V^{LC} \otimes (V^{RC})^{-}, W^{LC}\otimes (W^{RC})^{-}}$ and
a morphism $\Psi$ of $\C^{\times}$-rescalable partial pseudo-operads {}from 
$\tilde{K}^{c}\otimes \overline{\tilde{K}^{\bar{c}}}$ to the endomorphism 
partial pseudo-operad $H^{\C^{\times}}_{V^{LC}\otimes (V^{RC})^{-}, 
W^{LC}\otimes (V^{RC})^{-}}$.

\end{enumerate}

These data satisfy the following conditions:

\begin{enumerate}

\item For $s\in \R$ sufficiently negative, $V^{O}_{(s)}=0$ and 
for $m\in \Z$ sufficiently negative, $V^{LC}_{(m)}=V^{RC}_{(m)}=0$. 

\item For any $n\in \N$, $\Phi_{n}: \tilde{\mathfrak{S}}^{c}(n)\to 
H^{\R_{+}}_{V^{O}, W^{O}; V^{LC} \otimes V^{RC}, W^{LC}\otimes W^{RC}}(n)$ 
is linear on the fibers of 
$\tilde{\mathfrak{S}}^{c}(n)$.

\item The morphism $\Psi$ is equal to $\Psi^{L}\otimes \overline{\Psi^{R}}$
where $\Psi^{L}$ ($\Psi^{R}$) is a morphism of
$\C^{\times}$-rescalable partial pseudo-operads {}from 
$\tilde{K}^{c}$ ($\tilde{K}^{\bar{c}}$) to $H^{\C^{\times}}_{V^{LC},
W^{LC}}$ ($H^{\C^{\times}}_{V^{RC},
W^{RC}}$) and $\overline{\Psi^{R}}$ is the complex conjugate
of $\Psi^{R}$.  In addition, the triples $(V^{LC}, W^{LC}, \Psi^{L})$ and 
$(V^{RC}, W^{RC}, \Psi^{R})$ are meromorphic
$\tilde{K}^{c}$-associative algebra and $\tilde{K}^{\bar{c}}$-associative 
algebra, respectively.

\item For any $n\in \N$, the map 
$$\Phi_{n}: \tilde{\mathfrak{S}}^{c}(n)\to 
H^{\R_{+}}_{V^{O}, W^{O}; V^{LC} \otimes (V^{RC})^{-}, W^{LC}\otimes 
(W^{RC})^{-}}(n)$$
intertwines the action of the partial operad $\tilde{K}^{c}\otimes 
\overline{\tilde{K}^{\bar{c}}}$
on $\tilde{\mathfrak{S}}^{c}(n)$ and the action 
of the partial pseudo-operad $H^{\C^{\times}}_{V^{LC}\otimes (V^{RC})^{-}, 
W^{LC}\otimes (V^{RC})^{-}}$
on $H^{\R_{+}}_{V^{O}, W^{O}; 
V^{LC} \otimes (V^{RC})^{-}, W^{LC}\otimes (V^{RC})^{-}}(n)$.

\item For any $s_{1}, \dots, s_{n}\in \R$, there exists a finite 
subset $R(s_{1}, \dots, s_{n})\subset \R$ such that 
the image of 
$\coprod_{s\in s_{1}+\Z}V^{O}_{(s)}\otimes \cdots
\otimes \coprod_{s\in s_{n}+\Z}V^{O}_{(s)}$ under 
$\Phi_{n}(\psi^{\mathfrak{S}}_{n}(Q))$ for any $Q\in 
\tilde{\Upsilon}^{c}(n; 0)$
is in $\coprod_{s\in R(s_{1}, \dots, s_{n})+\Z}V^{O}_{(s)}$.

\item For any $v'\in (V^{O})'$, $u_{1}, \dots, u_{n}\in V^{O}$,
$v^{L}_{1}\otimes \bar{v}^{R}_{1}, \dots, v^{L}_{l}\otimes \bar{v}^{R}_{l}
\in V^{LC}\otimes 
(V^{RC})^{-}$, 
$$\langle v',
\Phi_{n}(\psi^{\mathfrak{S}}_{n}(Q))(u_{1}\otimes  \cdots \otimes u_{n}\otimes
v^{L}_{1}\otimes \bar{v}^{R}_{1}\otimes  \cdots \otimes 
v^{L}_{l}\otimes \bar{v}^{R}_{l})\rangle$$ 
as a function of
\begin{eqnarray*}
Q&=&(r_{1}, \dots, r_{n-1}; A^{(0)}, (a_{0}^{(1)}, A^{(1)}),
\cdots, (a_{0}^{(n)}, A^{(n)});\nn
&&\quad\quad\quad\quad\quad\quad  z_{1}, \dots, z_{l};
(b_{0}^{(1)}, B^{(1)}),
\cdots, (b_{0}^{(l)}, B^{(l)});\nn
&&\quad\quad\quad\quad\quad\quad  \bar{z}_{1}, \dots, \bar{z}_{l};
(\bar{b}_{0}^{(1)}, \bar{B}^{(1)}),
\cdots, (\bar{b}_{0}^{(l)}, \bar{B}^{(l)}))\nn
&&\quad\quad\quad\quad\quad\quad\in \tilde{\Upsilon}^{c}(n; l)
\subset K(n+2l)
\end{eqnarray*}
is of the 
form
\begin{eqnarray*}
\lefteqn{\sum_{i=1}^{k}f_{i}(r_{1}, \dots, r_{n-1}; z_{1}, \dots, z_{l};
\bar{z}_{1}, \dots, \bar{z}_{l})
\cdot}\nn
&&\cdot g_{i}(A^{(0)}, (a_{0}^{(1)}, A^{(1)}),
\cdots, (a_{0}^{(n)}, A^{(n)}); (b_{0}^{(1)}, B^{(1)}),
\cdots, (b_{0}^{(l)}, B^{(l)});\nn
&&\quad\quad\quad\quad\quad\quad\quad\quad\quad\quad\quad
\quad\quad\quad\quad\quad
 (\bar{b}_{0}^{(1)}, \bar{B}^{(1)}),
\cdots, (\bar{b}_{0}^{(l)}, \bar{B}^{(l)}))
\end{eqnarray*}
where the functions 
$$f_{i}(r_{1}, \dots, r_{n-1}; z_{1}, \dots, z_{l};
\xi_{1}, \dots, \xi_{l})$$
for $i=1, \dots, k$
are continuous differentiable in $r_{1}, \dots, r_{n-1}$
and are meromorphic in $z_{1}, \dots, z_{l}$,
$\xi_{1}, \dots, \xi_{l}$ with
$z_{i}=0, \infty$ and $z_{i}=z_{k}$, $i<k$, $z_{i}=r_{j}$,  
$\xi_{i}=0, \infty$, $\xi_{i}=\xi_{k}$, $i<k$, 
$\xi_{i}=r_{j}$ and $z_{i}=\xi_{k}$
for $i, k=1, \dots, l$ and $j=1, \dots, n-1$ as the only possible poles,
and 
\begin{eqnarray*}
\lefteqn{g_{i}(A^{(0)}, (a_{0}^{(1)}, A^{(1)}),
\cdots, (a_{0}^{(n)}, A^{(n)}); (b_{0}^{(1)}, B^{(1)}),
\cdots, (b_{0}^{(l)}, B^{(l)});}\nn
&&\quad\quad\quad\quad\quad\quad\quad\quad\quad\quad\quad
\quad\quad\quad\quad\quad
 (d_{0}^{(1)}, D^{(1)}),
\cdots, (d_{0}^{(l)}, D^{(l)}))
\end{eqnarray*}
for 
$i=1, \dots, k$ are polynomials in $A^{(0)}, (a_{0}^{(1)})^{\pm 1}, 
A^{(1)},
\cdots, (a_{0}^{(n)})^{\pm 1}, A^{(n)}$, $(b_{0}^{(1)})^{\pm 1}, 
B^{(1)},
\cdots, (b_{0}^{(l)})^{\pm 1}, B^{(l)}$ and $(d_{0}^{(1)})^{\pm 1}, 
D^{(1)},
\cdots, (d_{0}^{(l)})^{\pm 1}, D^{(l)}$.

\end{enumerate}

{\it Morphisms} (respectively, {\it isomorphisms}) of 
such pseudo-algebras
over $\tilde{\mathfrak{S}}^{c}$ 
are morphisms (respectively, isomorphisms) of the
underlying pseudo-algebras  over $\tilde{\mathfrak{S}}^{c}$ preserving 
all the structures. }
\end{defn}

We denote the differentiable-meromorphic pseudo-algebra 
over $\tilde{\mathfrak{S}}^{c}$
just defined by 
$$(V^{O}, W^{O}, V^{LC}, W^{LC}, V^{RC}, W^{RC}, \Phi, \Psi)$$
or simply by $(V^{O}, V^{LC}, V^{RC})$.
For these pseudo-algebras, 
we also have the following result whose proof is 
the same as those of the corresponding result in \cite{H5} and 
for Proposition \ref{pseudo}:

\begin{prop}
Any differentiable-meromorphic pseudo-algebra
over $\tilde{\mathfrak{S}}^{c}$
$$(V^{O}, W^{O}, V^{LC}, W^{LC}, V^{RC}, W^{RC}, \Phi, \Psi)$$
is an algebra over $\tilde{\mathfrak{S}}^{c}$, that is, the image 
of $\tilde{\mathfrak{S}}^{c}$ under $\Phi$ is a partial operad
(the image 
of $\tilde{\mathfrak{S}}^{c}$ under $\Phi$
satisfies  the composition-associativity).
\end{prop}

Because of this result, we shall omit the word ``pseudo'' {}from now on.

Note that given a vertex operator algebra $V$, its complex conjugate space 
$V^{-}$ has a natural vertex operator algebra structure 
$(V^{-}, Y^{-}, \overline{\one}, \overline{\omega})$
of central charge $\bar{c}$. 
We have the following generalization of Theorem \ref{main1}:

\begin{thm}\label{main2}
The category of objects consisting of a
grading-restricted conformal open-string vertex algebra
of central charge $c\in \C$, two
vertex operator algebras, one of central charge $c$ and the other 
of central charge $\bar{c}$,
and  homomorphisms {}from 
the first vertex operator algebra and the complex conjugate 
of the second vertex operator algebra to the meromorphic center of 
the grading-restricted conformal open-string vertex algebra
is isomorphic to the category of  
differentiable-meromorphic algebras
over $\tilde{\mathfrak{S}}^{c}$.
\end{thm}
\pf 
The proof of this theorem is basically the same as the proof of
isomorphism theorem for the geometric and operadic formulation of 
vertex operator algebras in \cite{H5} and the proof of 
Theorem \ref{main1}. The main new ingredient is that we use those
spheres with tubes which are obtained by doubling disks with 
strips and tubes. Here we give only a sketch.
More details will be given in \cite{K}. 

Given a grading-restricted conformal open-string vertex algebra 
$V^{O}$ of central charge $c$,
vertex operator algebras $V^{LC}$ and $V^{RC}$
of central charge $c$ and $\bar{c}$, respectively, 
and homomorphisms $h^{L}: V^{LC}\to C_{0}(V^{O})$
and $h^{R}: (V^{RC})^{-}\to C_{0}(V^{O})$. 
Let $W^{O}$, $W^{LC}$ and $W^{RC}$ be the modules 
for the Virasoro algebra generated by $\one^{O}$,  $\one^{LC}$ and
$\one^{RC}$ (the vacuums for these algebras), respectively. 
By the isomorphism theorem in \cite{H5}, 
there are meromorphic $\tilde{K}^{c}$-associative algebra
$(V^{LC}, W^{LC}, \Psi^{L})$ and meromorphic
$\tilde{K}^{\bar{c}}$-associative algebra
$(V^{RC}, W^{RC}, \Psi^{R})$ 
constructed {}from the vertex operator algebras
$V^{LC}$ and $V^{RC}$, respectively. Let $\Psi=\Psi^{L}\otimes 
\overline{\Psi^{R}}$.

By Theorem \ref{main1}, there is a 
differentiable $\tilde{\Upsilon}^{c}$-associative algebra
$(V^{O}, W^{O}, \Phi^{\Upsilon})$ constructed {}from $V^{O}$. 
Note that $\tilde{\Upsilon}^{c}$ is actually a partial suboperad of 
$\tilde{\mathfrak{S}}^{c}$. So this differentiable 
$\tilde{\Upsilon}^{c}$-associative algebra gives us part of 
a differentiable-meromorphic 
algebra structure over $\tilde{\mathfrak{S}}^{c}$. In general,
the construction of the differentiable-meromorphic 
algebra over $\tilde{\mathfrak{S}}^{c}$
can be obtained using the meromorphic $\tilde{K}^{c}$-associative algebra
$(V^{LC}, W^{LC}, \Psi^{L})$ and 
meromorphic $\tilde{K}^{c}$-associative algebra
$(V^{RC}, W^{RC}, \Psi^{R})$, 
the differentiable 
$\tilde{\Upsilon}^{c}$-associative algebra $(V^{O}, W^{O}, \Phi^{\Upsilon})$
and the homomorphisms $h^{L}$ and $h^{R}$. The action of 
the $\tilde{K}^{c}$-associative algebra 
$(V^{LC}\otimes (V^{RC})^{-}, W^{LC}\otimes (W^{RC})^{-}, \Psi)$
is also obtained using the homomorphisms $h^{L}$ and $h^{R}$. 
Thus we have a functor {}from the category of objects of the form 
$(V^{O}, V^{LC}, V^{RC}, h^{L}, h^{R})$
to the category of  
differentiable-meromorphic pseudo-algebra
over $\tilde{\mathfrak{S}}^{c}$.

Now given a differentiable-meromorphic pseudo-algebra
over $\tilde{\mathfrak{S}}^{c}$, by the definition and the isomorphism theorem 
in \cite{H5}, we know that there are vertex operator algebra structures
of central charge $c$ and $\bar{c}$
on $V^{LC}$ and $V^{RC}$, respectively. In particular, we have vacuums 
$\one^{LC}\in V^{LC}$
and $\one^{RC}\in V^{RC}$. Since $\tilde{\Upsilon}^{c}$ 
is actually a partial suboperad of 
$\tilde{\mathfrak{S}}^{c}$, by Theorem \ref{main1}, we also obtain a 
grading-restricted conformal open-string vertex algebra structure
of central charge $c$ on $V^{O}$. For $z\in \HH$, we
consider an
element $\Omega(z)$ of $\Upsilon(0; 1)$ which is the conformal 
equivalence class containing the following disk with strips and 
tubes of type $(1, 0; 0, 1)$: The disk $\hat{\HH}$ 
with the boundary puncture $\infty$ and the interior puncture $z$ 
and with the standard local coordinates vanishing at these punctures.
Then 
\begin{eqnarray*}
&\Phi_{0}(\psi^{\mathfrak{S}}_{0}(\Omega(z)))\in  
H^{\R_{+}}_{V^{O}, W^{O}; V^{LC}\otimes (V^{RC})^{-}, 
W^{LC}\otimes (W^{RC})^{-}}(0; 1)&\nn
&\quad\quad\quad \quad\quad\quad \quad\quad\quad \quad\quad\quad 
\quad\quad\quad =\hom(V^{LC}\otimes (V^{RC})^{-}, \overline{V^{O}}).&
\end{eqnarray*}
By Condition 6 in Definition \ref{sw-ch-alg}, 
we know that $\Phi_{0}(\psi^{\mathfrak{S}}_{0}(\Omega(z))(v^{L}\otimes 
\overline{\one^{RC}})$
is meromorphic in $z$ with the only pole $z=\infty$. In particular,
$$\lim_{z\to 0}\Phi_{0}(\psi^{\mathfrak{S}}_{0}(\Omega(z)))(v^{L}\otimes 
\overline{\one^{RC}})$$ 
exists. We define 
$$h^{L}(v^{L})=\lim_{z\to 0}\Phi_{0}(\psi^{\mathfrak{S}}_{0}(\Omega(z))
(v^{L}\otimes \overline{\one^{RC}}).$$
Thus we obtain
a linear map $h^{L}$ {}from $V^{LC}$ to $\overline{V^{O}}$. It is easy to see
that the image of $h^{L}$ is in fact in $C_{0}(V^{O})$ and 
$h^{L}$ is a homomorphism 
{}from $V^{LC}$ to $C_{0}(V^{O})$. Similarly we can construct
$h^{R}$.
Now we have a functor {}from the category of  
differentiable-meromorphic pseudo-algebra
over $\tilde{\mathfrak{S}}^{c}$ to the category of objects of the form 
$(V^{O}, V^{LC}, V^{RC}, h^{L}, h^{R})$.

{}From the isomorphism theorem in \cite{H5}, Theorem \ref{main1}
and the construction of the two functors above, 
we see that these two functors are inverse to each other. 
\epfv

In particular, we have:

\begin{cor} 
Let $V$ be a grading-restricted conformal open-string vertex
algebra of central charge $c$.  Then $V$ gives a natural structure of an
algebra over $\tilde{\mathfrak{S}}^{c}$ in the sense that
$(V, C_{0}(V), C_{0}(V)^{-})$ has a natural structure of 
differentiable-meromorphic 
algebra over $\tilde{\mathfrak{S}}^{c}$.
\end{cor}
\pf
We have a grading-restricted conformal algebra $V$ and a
vertex operator algebra $C_{0}(V)$. Let $h^{L}, h^{R}: C_{0}(V)\to C_{0}(V)$ 
be the identity map. Then $h^{L}$ and $h^{R}$ are homomorphisms {}from $C_{0}(V)$
to the meromorphic center of $V$. Then Theorem \ref{main2} gives
$(V, C_{0}(V), (C_{0}(V))^{-})$ a natural structure 
of differentiable-meromorphic 
algebra over $\tilde{\mathfrak{S}}^{c}$.
\epf

\noindent {\small \sc Department of Mathematics, Rutgers University,
110 Frelinghuysen Rd., Piscataway, NJ 08854-8019}

\vskip 1em

\noindent {\em E-mail address}: yzhuang@math.rutgers.edu, 
lkong@math.rutgers.edu

\end{document}